\documentclass{amsproc}

\title{Optimal Transport on Discrete Domains}
\author{Justin Solomon}
\address{MIT Department of Electrical Engineering and Computer Science}
\curraddr{}
\email{jsolomon@mit.edu}
\thanks{}

\subjclass[2010]{Primary: 65K10; Secondary: 52C99, 49Q20, 49M29}
\date{}

\usepackage{amsmath}
\usepackage{amsfonts}
\usepackage{amsthm}
\usepackage{longtable}
\usepackage{hyperref}
\usepackage{stmaryrd}
\usepackage{xcolor}
\usepackage{graphicx}
\hypersetup{
    colorlinks,
    linkcolor={red!50!black},
    citecolor={blue!50!black},
    urlcolor={blue!80!black}
}
\usepackage{units}
\usepackage{enumitem}
\usepackage{bbm}

\allowdisplaybreaks[4]

\newcommand{\R}[0]{\mathbb{R}}
\newcommand{\simplex}[0]{\mathbb{S}}
\newcommand{\W}[0]{\mathcal{W}}

\newcommand{\st}[0]{\textrm{s.t.}}

\newcommand{\KL}[0]{{\textrm{KL}}}
\newcommand{\diag}[0]{\mathrm{diag}}
\newcommand{\1}[0]{\mathbbm 1}
\newcommand{\OT}[0]{\mathrm{OT}}

\renewcommand{\L}[0]{\mathcal L}

\newcommand{\Prob}[0]{\mathrm{Prob}}
\newcommand{\Lip}[0]{\mathrm{Lip}}
\newcommand{\Lag}[0]{\mathrm{Lag}}

\newtheorem{theorem}{Theorem}[section]

\theoremstyle{definition}
\newtheorem{definition}[theorem]{Definition}

\theoremstyle{remark}

\numberwithin{equation}{section}

\graphicspath{ {figures/} }

\begin{document}

\begin{abstract}
Inspired by the matching of supply to demand in logistical problems, the optimal transport (or Monge--Kantorovich) problem involves the matching of probability distributions defined over a geometric domain such as a surface or manifold. In its most obvious discretization, optimal transport becomes a large-scale linear program, which typically is infeasible to solve efficiently on triangle meshes, graphs, point clouds, and other domains encountered in graphics and machine learning. Recent breakthroughs in numerical optimal transport, however, enable scalability to orders-of-magnitude larger problems, solvable in a fraction of a second. Here, we discuss advances in numerical optimal transport that leverage understanding of both discrete and smooth aspects of the problem. State-of-the-art techniques in discrete optimal transport combine insight from partial differential equations (PDE) with convex analysis to reformulate, discretize, and optimize transportation problems. The end result is a set of theoretically-justified models suitable for domains with thousands or millions of vertices. Since numerical optimal transport is a relatively new discipline, special emphasis is placed on identifying and explaining open problems in need of mathematical insight and additional research.
\end{abstract}

\maketitle

\section{Introduction}\label{sec:intro}

Many tools from discrete differential geometry (DDG) were inspired by practical considerations in areas like computer graphics and vision.  Disciplines like these require fine-grained understanding of geometric structure and the relationships between different shapes---problems for which the toolbox from smooth geometry can provide substantial insight.  Indeed, a triumph of discrete differential geometry is its incorporation into a wide array of computational pipelines, affecting the way artists, engineers, and scientists approach problem-solving across geometry-adjacent disciplines.

A key but neglected consideration hampering adoption of ideas in DDG in fields like computer vision and machine learning, however, is \emph{resilience} to noise and uncertainty.  The view of the world provided by video cameras, depth sensors, and other equipment is extremely unreliable.  Shapes do not necessarily come to a computer as complete, manifold meshes but rather may be scattered clouds of points that represent e.g.\ only those features visible from a single position.  Similarly, it may be impossible to pinpoint a feature on a shape exactly; rather, we may receive only a fuzzy signal indicating where a point or feature of interest \emph{may} be located.  Such uncertainty only increases in high-dimensional statistical contexts, where the presence of geometric structure in a given dataset is itself not a given.  Rather than regarding this messiness as an ``implementation issue'' to be coped with by engineers adapting DDG to imperfect data, however, the challenge of developing principled yet noise-resilient discrete theories of shape motivates new frontiers in mathematical research.

Probabilistic language provides a natural means of formalizing notions of uncertainty in the geometry processing pipeline.  In place of representing a feature or shape directly, we might instead use a probability distribution to encode a rougher notion of shape. Unfortunately, this proposal throws both smooth and discrete constructions off their foundations:  We must return to the basics and redefine notions like distance, distortion, and curvature in a fashion that does not rely on knowing shape with infinite precision and confidence.  At the same time, we must prove that the classical case is recovered as uncertainty diminishes to zero.

The mathematical discipline of \emph{optimal transport} (OT) shows promise for making geometry work in the probabilistic regime.  In its most basic form, optimal transport provides a means of lifting distances between points on a domain to distances between probability distributions \emph{over} the domain.  The basic construction of OT is to interpret probability distributions as piles of sand; the distance between two such piles of sand is defined as the amount of work it takes to transform one pile into the other.  This intuitive construction gave rise to an alternative name for OT in the computational world:  The ``earth mover's distance'' (EMD)~\cite{rubner2000earth}.  Indeed, the basic approach in OT is so natural that it has been proposed and re-proposed in many forms and with many names, from OT to EMD, the Mallows distance~\cite{levina2001earth}, the Monge--Kantorovich problem~\cite{villani2003topics}, the Hitchcock--Koopmans transportation problem~\cite{hitchcock1941distribution,koopmans1941exchange}, the Wasserstein/Vaser\v{s}te\u{i}n distance~\cite{vaserstein1969markov,dobrushin1970definition}, and undoubtedly many others.

Many credit Gaspard Monge with first formalizing the optimal transport problem in 1781 \cite{monge1781memoire}.  Beyond its early history, modern understanding of optimal transport dates back only to the World War II era, through the Nobel Prize-winning work of Leonid Kantorovich~\cite{kantorovich1942translocation}.   Jumping forward several decades, while many branches of DDG are dedicated to making centuries-old constructions on smooth manifolds work in the discrete case, optimal transport has the distinction of continuing to be an active area of research in the mathematical community whose basic properties are still being discovered.  Indeed, the computational and theoretical literature in this area move in lock-step:  New theoretical constructions often are adapted by the computational community in a matter of months, and some key theoretical ideas in transport were inspired by computational considerations and constructions.

Here, we aim to provide some intuition about transport and its relevance to the discrete differential geometry world.  While a complete survey of work on OT or even just its computational aspects is worthy of a full textbook, here we focus on the narrower problem of how to ``make transport work'' on a discretized piece of geometry amenable to representation on a computer.  The primary aim is to highlight the challenges in transitioning from smooth to discrete, to illustrate some basic constructions that have been proposed recently for this task, and---most importantly---to expose the plethora of open questions remaining in the relatively young discipline of computational OT.  No-doubt incomplete references are provided to selected intriguing ideas in computational OT, each of which is worthy of far more detailed discussion.

\subsubsection*{Additional reference.} Those readers with limited experience in related disciplines may wish to begin by reading~\cite{solomon2018computational}, a shorter survey by the author on the same topic, intended for a generalist audience.

\subsubsection*{Disclaimer.} These notes are intended as a short, intuitive, and \emph{extremely} informal introduction.  Optimal transport is a popular topic in mathematical research, and interested readers should refer to surveys such as~\cite{villani2003topics,villani2008optimal} for more comprehensive discussion.  The recent text~\cite{santambrogio2015optimal} provides discussion targeted to the applied world.  A few recent surveys also are targeted to computational issues in optimal transport~\cite{levy2017notions,peyre2017computational}.

The author of this tutorial offers his sincere apology to those colleagues whose foundational work is undoubtedly yet accidentally omitted from this document.  A ``venti''-sized caffeinated beverage is humbly offered in exchange for those readers' forgiveness and understanding.

\section{Motivation:  From Probability to Discrete Geometry}

To motivate the construction of optimal transport in the context of geometry processing, we begin by considering the case of smooth probability distributions over the real numbers $\R$.  Here, the geometry is extremely simple, described by values $x\in\R$ equipped with the distance metric $d(x,y):=|x-y|$.  Then we expand to define the transport problem in more generality and state a few useful properties.

\subsection{The Transport Problem}\label{sec:transport_problem}

Define the space of probability measures over $\R$ as $\Prob(\R)$. Without delving into the formalities of measure theory, these are roughly the functions $\mu\in\Prob(\R)$ assigning probabilities to sets $S\subseteq\R$ such that $\mu(S)\geq0$ for all measurable $S$, $\mu(\R)=1$, and $\mu(\cup_{i=1}^k S_i)=\sum_{i=1}^k\mu(S_i)$ for disjoint sets $\{S_i\subseteq\R\}_{i=1}^k$.  If $\mu$ is absolutely continuous, then it admits a \emph{distribution function} $\rho(x):\R\rightarrow\R$ assigning a probability density to every point:
$$\mu(S)=\int_S \rho(x)\,dx.$$

\begin{figure}[t]\centering
{\small\def\svgwidth{.9\textwidth}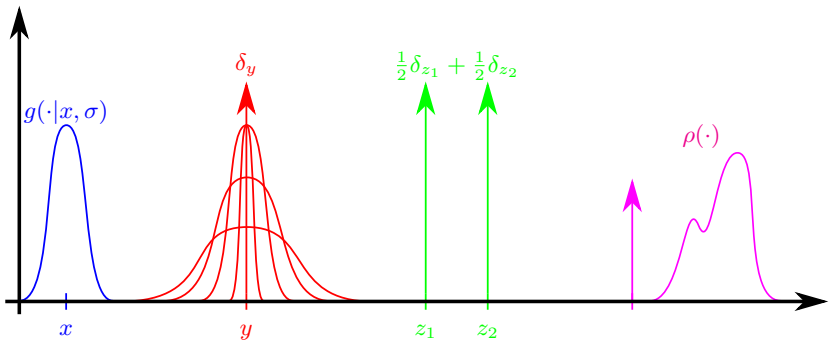}
\caption{One-dimensional examples of probability distributions used to encode geometric features with uncertainty. A probability distribution like a Gaussian {\color{blue}$g$} with standard deviation {\color{blue}$\sigma$} can be thought of as a ``fuzzy'' location of a point in {\color{blue}$x\in\R$}.  As a distribution sharpens about its mean to a $\delta$-function {\color{red}$\delta_y$}, it encodes a classical piece of geometry:  a point {\color{red}$y\in\R$}.  This language, however, is fundamentally broader, including constructions like the superposition of two points {\color{green}$z_1$} and {\color{green} $z_2$} or combining a point and a fuzzy feature into one distribution {\color{magenta}$\rho$}.}\label{fig:1d_example}
\end{figure}
Measure theory, probability, and statistics each are constructed from slightly different interpretations of the set of probability distributions $\Prob(\R)$.  Adding to the mix, we can think of optimal transport as a \emph{geometric} theory of probability.  In particular, as illustrated in Figure~\ref{fig:1d_example}, roughly a probability distribution over $\R$ can be thought of as a superposition of points in $\R$, whose weights are determined by $\rho(x)$.  We can recover a (complicated) representation for a single point $x\in\R$ as a Dirac $\delta$-measure centered at $x$.  

From a physical perspective, we can think of distributions geometrically using a physical analogy.  Suppose we are given a bucket of sand whose total mass is one pound.  We could distribute this pound of sand across the real numbers by stacking it all at a single point, concentrating it at a few points, or spreading it out smoothly.  The height of the pile of sand expresses a geometric feature:  Lots of sand at a point $x\in\R$ indicates we think a feature is located at $x$.

\begin{figure}\centering
{\small\def\svgwidth{.75\textwidth}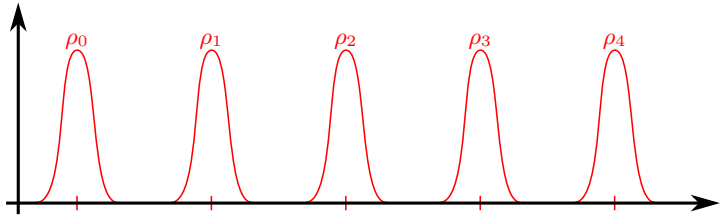}
\caption{The distributions {\color{red}$\rho_0,\ldots,\rho_4$} are equidistant with respect to the $L_1$ and KL divergences, while the Wasserstein distance from optimal transport increases linearly with distance over $\R$.}\label{fig:klbad}
\end{figure}

If we wish to deepen this analogy and lift notions from geometry to the space $\Prob(\R)$, perhaps the most basic object we must define is a notion of $\emph{distance}$ between two distributions $\mu_0,\mu_1\in\Prob(\R)$ that resembles the distance $d(x,y)=|x-y|$ between points on the underlying space.  Supposing for now that $\mu_0$ and $\mu_1$ admit distribution functions $\rho_0$ and $\rho_1$, respectively, a few candidate notions of distance or divergence come to mind:
\begin{align*}
\textrm{$L_1$ distance:} & \hspace{.25in}d_{L_1}(\rho_0,\rho_1):=\int_{-\infty}^\infty |\rho_0(x)-\rho_1(x)|\,dx\\
\textrm{KL divergence:} & \hspace{.25in}d_{\mathrm{KL}}(\rho_0\|\rho_1):=\int_{-\infty}^\infty \rho_0(x)\log\frac{\rho_0(x)}{\rho_1(x)}\,dx.
\end{align*}
These divergences are used widely in analysis and information theory, but they are insufficient for geometric computation.  In particular, consider the distributions in Figure~\ref{fig:klbad}.    The two divergences above give the same value for any pair of different $\rho_i$'s!  This is because they measure only the overlap; the ground distance $d(x,y)=|x-y|$ is never used in their computation.

Optimal transport resolves this issue by leveraging the physical analogy proposed above.  In particular, suppose our sand is currently in arrangement $\rho_0$ and we wish to \emph{reshape} it to a new distribution $\rho_1$.  We take a steam shovel and begin scooping up the sand at points $x$ in $\rho_0$ where $\rho_0(x)>\rho_1(x)$ and dropping it places where $\rho_1(x)>\rho_0(x)$; eventually one distribution is transformed into the other.

\begin{figure}\centering
\begin{tabular}{cc}
{\small\def\svgwidth{.3\textwidth}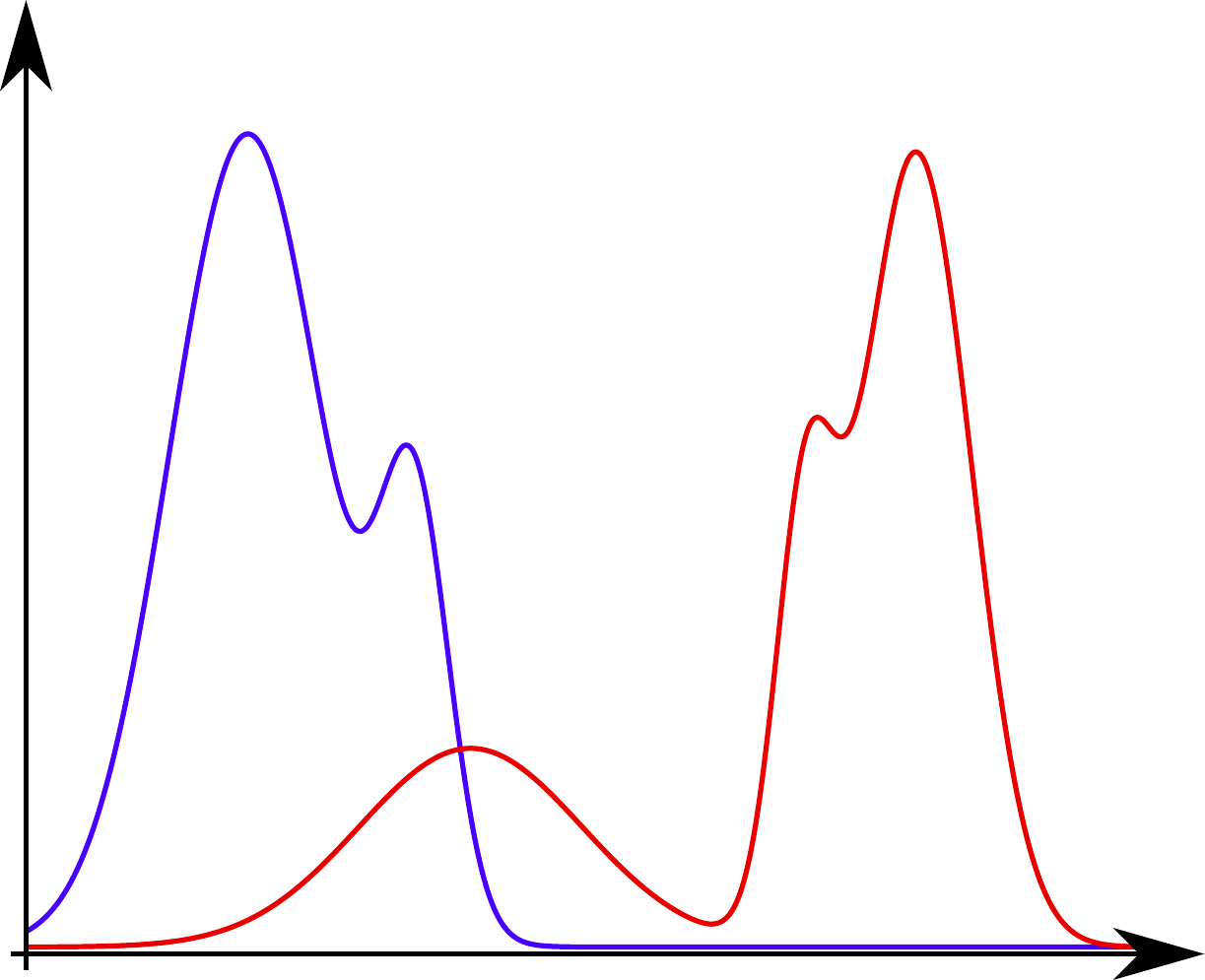}&
{\small\def\svgwidth{.3\textwidth}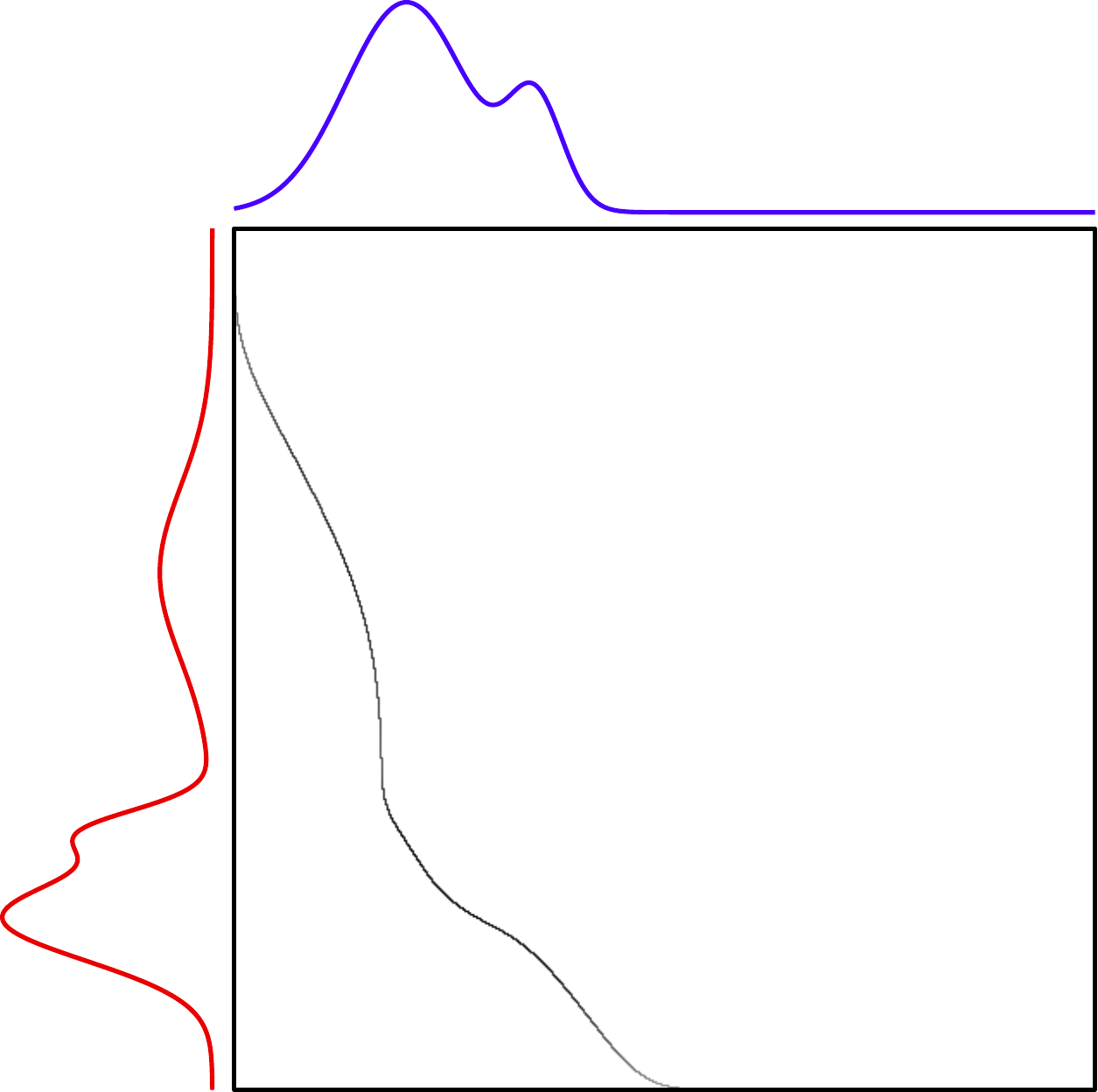}\\
(a) Source and target &\hspace{.1in} (b) Transport map
\end{tabular}
\caption{Two distributions over the real line (a) and the resulting transport map (b).  In (b), the box is the product space $[0,1]\times[0,1]$, and dark values indicate a matching between $\rho_0$ and $\rho_1$.}\label{fig:1dmap}
\end{figure}

There are many ways the steam shovel could approach its task:  We could move sand efficiently, or we could drive it miles away and then drive back, wasting fuel in the process.  But assuming $\rho_0\neq\rho_1$, there is some amount of work inherent in the fact that $\rho_0$ and $\rho_1$ are not the same.  We can formalize this idea by solving for an unknown measure $\pi(x,y)$ determining how much mass gets moved from $x$ to $y$ by the steam shovel for each $(x,y)$ pair.  The minimum amount of work is then
\begin{equation}\label{eq:w1_1d}
\W_1(\rho_0,\rho_1):=\left\{\!
\begin{array}{r@{\ }ll}
\min_\pi & \iint_{\R\times\R} \pi(x,y)|x-y|\,dx\,dy & \textrm{Minimize total work}\\
\st & \pi\geq0\ \forall x,y\in \R & \textrm{Nonnegative mass}\\
& \int_\R \pi(x,y)\,dy = \rho_0(x)\,\forall x\in\R & \textrm{Starts from $\rho_0$}\\
& \int_\R \pi(x,y)\,dx = \rho_1(y)\,\forall y\in\R & \textrm{Ends at $\rho_1$}.
\end{array}\right.
\end{equation}
This optimization problem quantifies the minimum amount of work---measured as mass $\pi(x,y)$ times distance traveled $|x-y|$---required to transform $\rho_0$ into $\rho_1$.  We can think of the unknown function $\pi$ as the instructions given to the laziest possible steam shovel tasked with dropping one distribution onto another.  This amount of work is known as the \emph{1-Wasserstein distance} in optimal transport; in one dimension, it equals the $L_1$ distance between the cumulative distribution functions of $\rho_0$ and $\rho_1$.  An example of $\rho_0$, $\rho_1$, and the resulting $\pi$ is shown in Figure~\ref{fig:1dmap}.

Generalizing slightly, we can define the $p$-Wasserstein distance:
\begin{equation}\label{eq:wp_1d}
[\W_p(\rho_0,\rho_1)]^p:=\left\{
\begin{array}{rl}
\min_\pi & \iint_{\R\times\R} \pi(x,y)|x-y|^p\,dx\,dy \\
\st & \pi\geq0\ \forall x,y\in \R \\
& \int_\R \pi(x,y)\,dy = \rho_0(x)\,\forall x\in\R \\
& \int_\R \pi(x,y)\,dx = \rho_1(y)\,\forall y\in\R .
\end{array}\right.
\end{equation}
In analogy to Euclidean space, many properties of $\W_p$ are split into cases $p<1$, $p=1$, and $p>1$; for instance, it satisfies the triangle inequality any time $p\geq1$. The $p=2$ case is of particular interest in the literature and corresponds to a ``least-squares'' version of transport that minimizes kinetic energy rather than work (see \S\ref{sec:manyformulas}).  Generalizing~\eqref{eq:wp_1d} even more, if we replace $|x-y|^p$ with a generic cost $c(x,y)$ we recover the \emph{Kantorovich problem}~\cite{kantorovich1942translocation}.

It is important to note an alternative formulation of the transport problem~\eqref{eq:wp_1d}, which historically was posed first but does not always admit a solution.  Rather than optimizing for a function $\pi(x,y)$ with an unknown for every possible $(x,y)$ pair, one could consider an alternative in which instead the variable is a single function $\phi(x)$ that ``pushes forward'' $\rho_0$ onto $\rho_1$; this corresponds to choosing a single destination $\phi(x)$ for every source point $x$.  In this case, the objective function would look like
\begin{equation}\label{eq:monge}
\int_{-\infty}^\infty |\phi(x)-x|^p\rho_0(x)\,dx,
\end{equation}
and the constraints would ask that the image of $\rho_0$ under $\phi$ is $\rho_1$, notated $\phi_\sharp\rho_0=\rho_1$. While this version corresponds to the original version of transport proposed by Monge, sometimes for the transport problem to be solvable it is necessary to split the mass at a single source point to multiple destinations.  A triumph of theoretical optimal transport, however, shows that $\pi(x,y)$ is nonzero only on some set $\{(x,\phi(x)):x\in\R\}$ whenever $\rho_0$ is absolutely continuous, linking the two problems.

\subsection{Discrete Problems in One Dimension}
% discrete/semidiscrete/continuous

So far our definitions have not been amenable to numerical computation:  Our unknowns are functions $\pi(x,y)$ with \emph{infinite} numbers of variables (one value of $\pi$ for each $(x,y)$ pair in $\R\times\R$)---certainly more than can be stored on a computer with finite capacity.  Continuing to work in one dimension, we suggest some special cases where we can solve the transport problem with a finite number of variables.

Rather than working with distribution functions $\rho(x)$, we will relax to the more general case of transport between measures $\mu_0,\mu_1\in\Prob(\R)$.  Define the Dirac $\delta$-measure centered at $x\in\R$ via
$$
\delta_x(S):=\left\{
\begin{array}{rl}
1 & \textrm{ if }x\in S\\
0 & \textrm{ otherwise.}
\end{array}
\right.
$$
It is easy to check that $\delta_x(\cdot)$ is a probability measure.

\begin{figure}\centering
\begin{tabular}{cc}
{\small\def\svgwidth{.3\textwidth}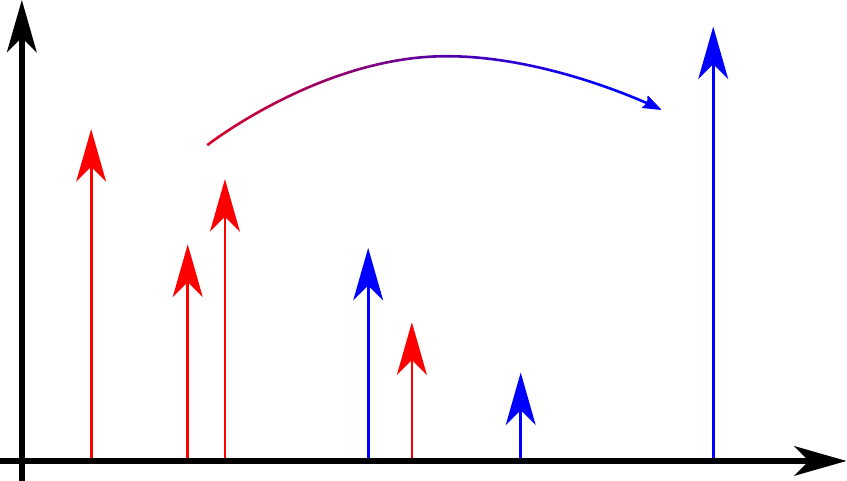}&
{\small\def\svgwidth{.3\textwidth}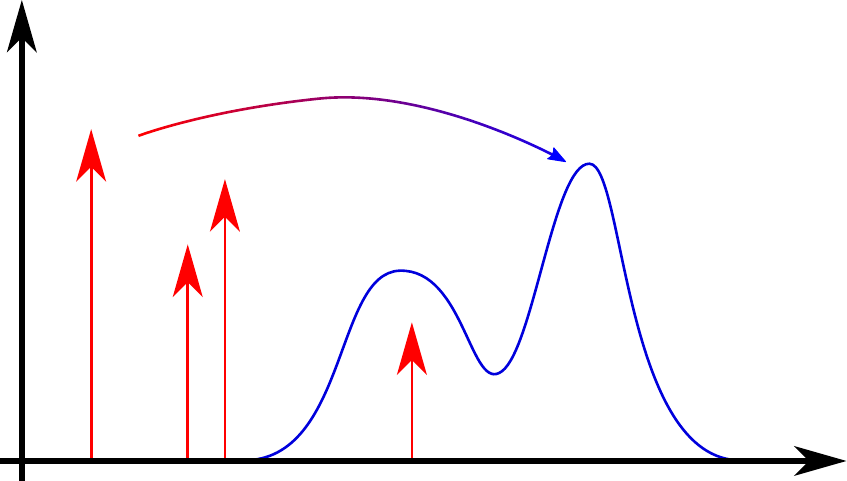}\\
(a) Fully discrete transport &\hspace{.1in} (b) Semidiscrete transport
\end{tabular}
\caption{Discrete (a) and semidiscrete (b) optimal transport in one dimension.}\label{fig:transportexamples}
\end{figure}

Suppose $\mu_0,\mu_1\in\Prob(\R)$ can be written as \emph{superpositions} of $\delta$ measures:
\begin{equation}\label{eq:deltasum}
\mu_0:=\sum_{i=1}^{k_0} a_{0i} \delta_{x_{0i}}\hspace{.5in}\textrm{and}\hspace{.5in}
\mu_1:=\sum_{i=1}^{k_1} a_{1i} \delta_{x_{1i}},
\end{equation}
where $1=\sum_{i=1}^{k_0}a_{0i}=\sum_{i=1}^{k_1}a_{1i}$ and $a_{0i},a_{1i}\geq0$ for all $i$.  Figure~\ref{fig:transportexamples}(a) illustrates this case; all the mass of $\mu_0$ and $\mu_1$ is concentrated at a few isolated points.

In the case where the source and target distributions are composed of $\delta$'s, we only can move mass between pairs of points $x_{0i}\mapsto x_{1j}$.  Taking $T_{ij}$ the total mass moved from $x_{0i}$ to $x_{1j}$, we can solve for $\W_p^p$ as
\begin{equation}\label{eq:wp_finite}
[\W_p(\mu_0,\mu_1)]^p=\left\{
\begin{array}{rl}
\min_{T\in\R^{k_0\times k_1}} & \sum_{ij} T_{ij} |x_{0i}-x_{1j}|^p\\
\st & T\geq0\\
& \sum_j T_{ij} = a_{0i}\\
& \sum_i T_{ij} = a_{1j}.
\end{array}
\right.
\end{equation}
This is an optimization problem in $k_0k_1$ variables $T_{ij}$:  No need for an infinite number of $\pi(x,y)$'s!  In fact, it is a \emph{linear program} solvable using many classic algorithms, such as the simplex or interior point methods.%  In fact, this provides a \emph{convergent} means of approximating transport distances, in the sense that any measure can be well-approximated by a sum of the form~\eqref{eq:deltasum} in the Wasserstein metric~\cite{kloeckner2012approximation}.

There is a more subtle case where we can still represent the unknown in optimal transport using a finite number of variables.  Suppose $\mu_0\in\Prob(\R)$ is a superposition of $\delta$ measures and $\mu_1\in\Prob(\R)$ is absolutely continuous, implying $\mu_1$ admits a distribution function $\rho_1(x)$:
\begin{equation}\label{eq:deltasum2}
\mu_0:=\sum_{i=1}^{k} a_{i} \delta_{x_{i}}\hspace{.5in}\textrm{and}\hspace{.5in}
\mu_1(S):=\int_S \rho_1(x)\,dx.
\end{equation}
This situation is illustrated in Figure~\ref{fig:transportexamples}(b); it corresponds to transporting from a distribution whose mass is concentrated at a few points to a distribution whose distribution is more smooth.  In the technical literature, this setup is known as \emph{semidiscrete} transport. 

Returning to the transport problem in~\eqref{eq:wp_1d}, in this semidiscrete case we can think of the coupling $\pi$ as decomposing into a set of measures $\pi_1,\pi_2,\ldots,\pi_k\in\Prob(\R)$ where each term in the sum~\eqref{eq:deltasum2} has its own target distribution:  $\delta_{x_i}\mapsto \pi_i.$  As a sanity check, note that $\mu_1=\sum_i a_i\pi_i(x).$

\begin{figure}\centering
{\small\def\svgwidth{.5\textwidth}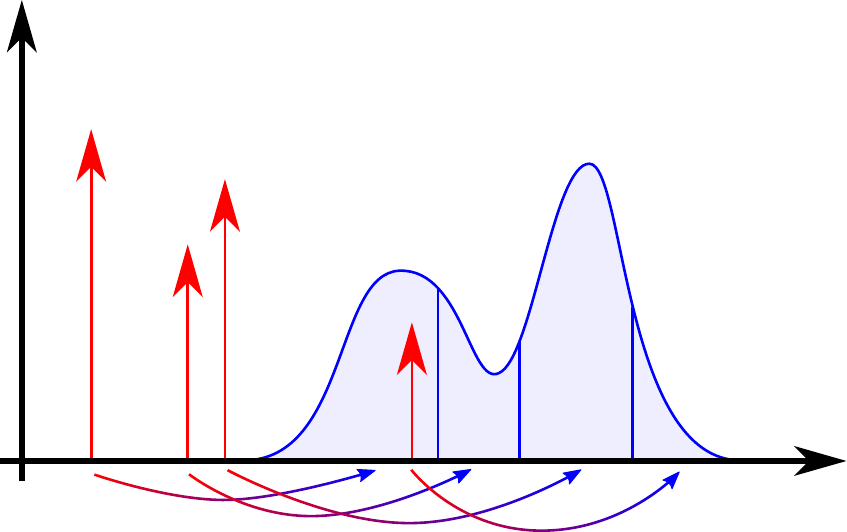}
\caption{Solving 1D semidiscrete transport from Figure~\ref{fig:transportexamples}(b); every Dirac $\delta$-function mass in the source {\color{red}$\mu_0$} gets mapped to a contiguous interval worth of mass in the target {\color{blue}$\mu_1$}.}\label{fig:1dsemidiscretesolution}
\end{figure}

Without loss of generality, we can assume the $x_i$'s are sorted, that is, $x_1<x_2<\cdots<x_k$.  Suppose $1\leq i<j\leq k$, and hence $x_i<x_j$.  In one dimension, it is easy to see that the optimal transport map $\pi$ should never ``leapfrog'' mass, that is, the delivery target of the mass at $x_i$ when transported to $\rho_1$ should be \emph{to the left} of the target of mass at $x_j$, as illustrated in Figure~\ref{fig:1dsemidiscretesolution}.  This monotonicity property implies the existence of intervals $[b_1,c_1],[b_2,c_2],\ldots,[b_k,c_k]$ such that $\pi_i$ is supported in $[b_i,c_i]$ and $c_i<b_j$ whenever $i<j$; the mass $a_i\delta_{x_i}$ is distributed according to $\rho_1(x)$ in the interval $[b_i,c_i]$.  

The semidiscrete transport problem corresponds to another case where we can solve a transport problem with a finite number of variables, the $b_i$'s and $c_i$'s.  Of course, in one dimension these can be read off from the cumulative distribution function (CDF) of $\rho_1$, but in higher dimensions this will not be the case.  Instead, the intervals $[b_i,c_i]$ will be replaced with \emph{power cells}, a generalization of a Voronoi diagram (\S\ref{sec:semidiscrete_tr}).

While our discussion above gives two cases in which a computer could plausibly solve the transport problem, they do not correspond to the usual situation for DDG in which the geometry itself---in this case the real line $\R$---is discretized.  As we will see in the discussion in future sections, there currently does not exist consensus about what to do in this case but several possible adaptations to this case have been proposed.

\subsection{Moving to Higher Dimensions}

We are now ready to state the optimal transport problem in full generality.  Following~\cite[\S1.1.1]{villani2003topics}, take $(X,\mu)$ and $(Y,\nu)$ to be probability spaces, paired with a nonnegative measurable cost function $c(x,y)$.  Define a \emph{measure coupling} $\pi\in\Pi(\mu,\nu)$ as follows:
\begin{definition}[Measure coupling]
A \emph{measure coupling} $\pi\in\Prob(X\times Y)$ is a probability measure on $X\times Y$ satisfying
\begin{align*}
\pi(A\times Y)&=\mu(A)\\
\pi(X\times B)&=\nu(B)
\end{align*}
for all measurable $A\subseteq X$ and $B\subseteq Y$.  The set of measure couplings between $\mu$ and $\nu$ is denoted $\Pi(\mu,\nu)$.
\end{definition}

With this piece of notation, we can write the Kantorovich optimal transport problem as follows:
\begin{equation}\label{eq:transport}
\boxed{
\OT(\mu,\nu;c):=
%\left\{
%\begin{array}{rl}
%\min_{\pi\in\Prob(X\times Y)} & \iint_{X\times Y}c(x,y)\,d\pi(x,y)\\
%\st & \pi\in\Pi(\mu,\nu).
%\end{array}\right.
\min_{\pi\in\Pi(\mu,\nu)}\ \iint_{X\times Y}c(x,y)\,d\pi(x,y)
}
\end{equation}
Here, we use some notation from measure theory:  $d\pi(x,y)$ denotes integration against probability measure $\pi$.  Note if $\pi$ admits a distribution function $p(x,y)$ then we can write $d\pi(x,y)=p(x,y)\,dx\,dy$; the more general  notation allows for $\delta$ measures and other objects that cannot be written as functions.

We note a few interesting special cases below:

\subsubsection*{Discrete transportation.}  Suppose $X=\{1,2,\ldots,k_1\}$ and $Y=\{1,2,\ldots,k_2\}$.  Then, $\mu\in\Prob(X)$ can be written as a vector $v\in\simplex_{k_1}$ and $\nu\in\Prob(Y)$ can be written as a vector $w\in\simplex_{k_2}$, where $\simplex_k$ denotes the $k$-dimensional probability simplex:
\begin{equation}
\simplex_k:=\left\{v\in\R^k : v\geq0\textrm{ and }\sum_i v_i=1\right\}.
\end{equation}
Our cost function becomes discrete as well and can be written as a matrix $C=(c_{ij})$.  After simplification, the transport problem between $v\in\simplex_{k_1}$ and $w\in\simplex_{k_2}$ given cost matrix $C$ becomes
\begin{equation}\label{eq:discrete_ot}
\OT(v,w;C)=\left\{
\begin{array}{rl}
\min_{T\in\R^{k_1\times k_2}} & \sum_{ij} T_{ij}c_{ij}\\
\st & T\geq0\\
& \sum_j T_{ij}=v_i\ \forall i\in\{1,\ldots,k_1\}\\
& \sum_i T_{ij}=w_j\ \forall j\in\{1,\ldots,k_2\}.
\end{array}
\right.
\end{equation}
This linear program is solvable computationally and is the most obvious way to make optimal transport work in a discrete context.  It was proposed in the computational literature as the ``earth mover's distance'' (EMD)~\cite{rubner2000earth}.  When $k_1=k_2:=k$ and $C$ is symmetric, nonnegative, and satisfies the triangle inequality, one can check that $\OT(\cdot,\cdot;C)$ is a distance on $\simplex_k$; see~\cite{cuturi2014ground} for a clear proof of this property.

\subsubsection*{Wasserstein distance.}  Next, suppose $X=Y=\R^n$, and take $c_{n,p}(x,y):=\|x-y\|_2^p$.  Then, we recover the \emph{Wasserstein distance} on $\Prob(\R^n)$, defined via
\begin{equation}\label{eq:wasserstein}
\W_p(\mu,\nu):=[\OT(\mu,\nu;c_{n,p})]^{\nicefrac{1}{p}}.
\end{equation}
$\W_p$ is a distance when $p\geq1$, and $\W_p^p$ is a distance when $p\in[0,1)$~\cite[\S7.1.1]{villani2003topics}.  In fact, the Wasserstein distance can be defined for probability measures over a surface, Riemannian manifold, or even a Polish space via the same formula.

The Wasserstein distance has drawn considerable application-oriented interest and aligns well with the basic motivation laid out in \S\ref{sec:intro}.  Its basic role is to lift distances between points $\|x-y\|_2^p$ to distances between distributions in a compatible fashion:  The Wasserstein distance between two $\delta$-functions $\delta_x$ and $\delta_y$ is exactly the distance $\|x-y\|_2$.  In \S\ref{sec:app}, we will show how this basic property has strong bearing on several computational pipelines that need to lift geometric constructions to uncertain contexts.

\subsection{One Value, Many Formulas}\label{sec:manyformulas}

A remarkable property of the transport problem~\eqref{eq:transport} is the sheer number of equivalent formulations that all lead to the same value, the cost of transporting mass from one measure onto another.  These not only provide many interpretations of the transport problem but also suggest a diverse set of computational algorithms for transport, each of which tackles a different way of writing down the basic problem.

\subsubsection*{Duality.}  A basic idea in the world of convex optimization is that of \emph{duality}, that every minimization problem admits a ``dual'' maximization problem whose optimal value lower-bounds that of the primal.  As with most linear programs, optimal transport typically exhibits \emph{strong duality}:  The optimal values of the maximization and minimization problems coincide.

To motivate duality for transport, we will start with the finite-dimensional problem~\eqref{eq:discrete_ot}.  We note two simple identities:
$$
\max_{s\in\R} st = \left\{
\begin{array}{cl}
0 & \textrm{ if }t=0\\
\infty & \textrm{ otherwise}
\end{array}
\right.
\hspace{1in}
\max_{s\leq0} st = \left\{
\begin{array}{cl}
0 & \textrm{ if }t\geq0\\
\infty & \textrm{ otherwise}
\end{array}
\right.
$$
These allow us to write~\eqref{eq:discrete_ot} as follows:
$${\color{red}\min_T} {\color{blue}\max_{S\leq0,\phi,\psi}} \left[\sum_{ij} T_{ij}(c_{ij}+S_{ij})+\sum_i \phi_i \left(v_i-\sum_jT_{ij}\right) + \sum_j \psi_j\left(w_j-\sum_i T_{ij}\right)\right].$$
The dual problem is derived by simply swapping the min and the max:
\begin{align*}{\color{blue}\max_{S\leq0,\phi,\psi}} {\color{red}\min_T} &\left[\sum_{ij} T_{ij}(c_{ij}+S_{ij})+\sum_i \phi_i \left(v_i-\sum_jT_{ij}\right) + \sum_j \psi_j\left(w_j-\sum_i T_{ij}\right)\right]\\
&=\max_{S\leq0,\phi,\psi}\min_T \left[\sum_{ij} T_{ij}(c_{ij}+S_{ij}-\phi_i-\psi_j)+\sum_i \phi_i v_i + \sum_j \psi_jw_j\right]\textrm{ after refactoring.}
\end{align*}
Since $T$ is unbounded in the inner optimization problem of the dual, the solution of the inner minimization is $-\infty$ unless $S_{ij}=\phi_i+\psi_j-c_{ij}$ for all $(i,j)$, that is, unless the coefficient of $T_{ij}$ equals zero.  Since the outer problem is a maximization, clearly we should avoid an optimal value of $-\infty$ for the inner minimization.  Hence, we can safely add $S_{ij}=\phi_i+\psi_j-c_{ij}$ as a constraint to the dual problem.  After some simplification, we arrive at the dual of~\eqref{eq:discrete_ot}:
\begin{equation}\label{eq:discrete_dual}
\begin{array}{rl}
\max_{\phi,\psi} & \sum_i [\phi_i v_i + \psi_i w_i]\\
\st & \phi_i + \psi_j \leq c_{ij} \ \forall (i,j).
\end{array}
\end{equation}
Although we have not justified that it is acceptable to swap a max and a min in this context, several techniques ranging from direct proof to the ``sledgehammer'' Slater duality condition~\cite{slater1950lagrange} show that the optimal value of this maximization problem agrees with the optimal value of the minimization problem~\eqref{eq:discrete_ot}.

As is often the case in convex optimization, the dual~\eqref{eq:discrete_dual} of the transport problem~\eqref{eq:discrete_ot} has an intuitive interpretation.  Suppose we change roles in optimal transport from the worker who wishes to minimize work to a company that wishes to maximize profit.  The customer pays $\phi_i$ dollars per pound to drop off material $v_i$ to ship from location $i$ and $\psi_j$ dollars per pound to pick up material $w_j$ from location $j$.  The dual problem~\eqref{eq:discrete_dual} maximizes profit under the constraint that it is never cheaper for the customer to just drive from $i$ to $j$ and ignore the service completely:  $\phi_i+\psi_j\leq c_{ij}.$

We pause here to note some rough trade-offs between the primal and dual transport problems.  Since both yield the same optimal value, the designer of a computational system for solving optimal transport problems has a decision to make:  whether to solve the primal problem, the dual problem, or both simultaneously (the latter aptly named a ``primal--dual'' algorithm).  There are advantages and disadvantages to each approach.  The primal problem~\eqref{eq:discrete_ot} directly yields the matrix $T$, which tells not just the cost of transport but how much mass $T_{ij}$ to move from source $i$ to destination $j$; the only inequality constraint is that the entire matrix has nonnegative entries.  On the other hand, the dual problem~\eqref{eq:discrete_dual} has fewer variables, making it easier to store the output on the computer, but the ``shadow price'' variables $(\phi,\psi)$ are harder to interpret and are constrained by a quadratic number of inequalities.  Currently there is little consensus as to which formulation leads to more successful algorithms or discretizations, and state-of-the-art techniques are divided among the two basic approaches.

As with many constructions in optimal transport, the dual of the measure-theoretic problem~\eqref{eq:transport} resembles the discrete case up to a change of the notation.  In particular, we can write 
\begin{equation}\label{eq:transportdual}
\boxed{
\OT(\mu,\nu;c):=\left\{
\begin{array}{rl}
\sup_{\substack{\phi\in L^1(d\mu),\\\psi\in L^1(d\nu)}} & \int_X \phi(x)\,d\mu(x) + \int_Y \psi(y)\,d\nu(y)\\
\st & \phi(x)+\psi(y)\leq c(x,y)\\ &\ \ \ \ \textrm{ for }d\mu\textrm{-a.e.\ }x\in X,\ d\nu\textrm{-a.e.\ }y\in Y.
\end{array}\right.
}
\end{equation}

It is worth noting a simplification that appears often in the transport world.  Since $\mu$ and $\nu$ are positive measures and the overall problem in~\eqref{eq:transportdual} is a maximization, we might as well choose $\phi$ and $\psi$ as large as possible while satisfying the constraints.  Suppose we fix the function $\phi(x)$ and \emph{just} optimize for the function $\psi(x)$.  Rearranging the constraint shows that for all $(x,y)\in X\times Y$ we must have $\psi(y)\leq c(x,y)-\phi(x).$  Equivalently, for all $y\in Y$ we must have $\psi(y)\leq\inf_{x\in X} [c(x,y)-\phi(x)]$.  Define the \emph{$c$-transform} 
\begin{equation}\label{eq:c_transform}
\phi^c(y):=\inf_{x\in X} [c(x,y)-\phi(x)].
\end{equation}
By the argument above we have
$$
\OT(\mu,\nu;c)=
\sup_{\phi\in L^1(d\mu)}  \int_X \phi(x)\,d\mu(x) + \int_Y \phi^c(y)\,d\nu(y).
$$
This problem is unconstrained, but the transformation from $\phi$ to $\phi^c$ is relatively complicated.

We finally note one special case of this dual formula, the 1-Wasserstein distance, which has gained recent interest in the machine learning world thanks to its application in generative adversarial networks (GANs)~\cite{arjovsky2017wasserstein}.  In this case, $X=Y=\R^n$ and $c(x,y)=\|x-y\|_2$.  We can derive a bound as follows:
\begin{align}
|\phi^c(x)-\phi^c(y)|
&=\left|
\inf_{z} [\|x-z\|_2-\phi(z)] - \inf_{z} [\|y-z\|_2 - \phi(z)]
\right|\textrm{ by definition}\nonumber\\
&\leq\sup_z |\|x-z\|_2-\|y-z\|_2|\nonumber\\&\hspace{.25in}\textrm{ by the identity }|\inf_x f(x)-\inf_x g(x)|\leq\sup_x |f(x)-g(x)|\nonumber\\
&\leq \|x-y\|_2\textrm{ by the reverse triangle inequality.}\label{eq:tri}
\end{align}
Furthermore, by definition of the $c$-transform~\eqref{eq:c_transform} by taking $x=y$ we have $\phi^c(y)\leq -\phi(y)$, or equivalently $\phi(y)\leq-\phi^c(y)$. Hence,
\begin{align*}
\W_1(\mu,\nu)
&=\OT(\mu,\nu;c)\textrm{ through our choice }c(x,y):=\|x-y\|_2\\
&=\sup_{\phi\in L^1(d\mu)}\int_{\R^n} \phi(x)\,d\mu(x)+\int_{\R^n}\phi^c(y)\,d\nu(y)\textrm{ by definition of the $c$-transform}\\
&\leq \int_{\R^n} \phi^c(x)\,[d\nu(x)-d\mu(x)]\textrm{ since }\phi(y)\leq-\phi^c(y)\ \forall y\in\R^n\\
&\leq \sup_{\psi\in\Lip_1(\R^n)}\int_{\R^n}\psi(x)\,[d\nu(x)-d\mu(x)]\\
&\hspace{.2in}\textrm{ where }\Lip_1(\R^n):=\{f(x): |f(x)-f(y)|\leq\|x-y\|_2\ \forall x,y\in\R^n\}.
\end{align*}
$\Lip_1$ denotes the set of $1$-Lipschitz functions; the last step is derived from~\eqref{eq:tri}, which shows that $\psi^c$ is $1$-Lipschitz.

In fact, this inequality is an equality.  To prove this, take $\psi$ to be any 1-Lipschitz function.  Then, 
\begin{equation}\label{eq:lipbound}
\psi^c(y)=\inf_{x\in\R^n} [\|x-y\|_2-\psi(x)]\geq\inf_{x\in\R^n} [\|x-y\|_2-\|x-y\|_2-\psi(y)]=-\psi(y).
\end{equation}
where we have rearranged the Lipschitz property $\psi(x)-\psi(y)\leq\|x-y\|_2$ to show $-\psi(x)\geq-\|x-y\|_2-\psi(y)$.  Hence,  
\begin{align*}
\sup_{\psi\in\Lip_1(\R^n)}\int_{\R^n}\hspace{-.1in}\psi(x)\,[d\nu(x)\!-\!d\mu(x)]
&\leq \sup_{\psi\in\Lip_1(\!\R^n\!)}\int_{\R^n}\hspace{-.1in}[\psi(x)\,d\nu(x)+\psi^c(y)]\, d\mu(y)\textrm{ by~\eqref{eq:lipbound}}\\
&\leq \sup_{\psi\in L^1(\!d\nu\!)}\int_{\R^n}\hspace{-.1in}[\psi(x)\,d\nu(x)+\psi^c(y)]\, d\mu(y)\\&\hspace{.25in}\textrm{ since the constraints are loosened}\\
&=\W_1(\mu,\nu).
\end{align*}
This finishes motivating our final formula
$$
\boxed{
\W_1(\mu,\nu)=\sup_{\psi\in\Lip_1(\R^n)}\int_{\R^n}\psi(x)\,[d\nu(x)-d\mu(x)].
}
$$
This convenient identity is used in computational contexts because the constraint that a function is 1-Lipschitz is fairly easy to enforce computationally; sadly, it does not extend to other Wasserstein $\W_p$ distances, which have nicer uniqueness and regularity properties when $p>1$.

\subsubsection*{Eulerian transport.}  The language of fluid dynamics introduces two equivalent ways to understand the flow of a liquid or gas as it sloshes in a tank.  In the \emph{Lagrangian} framework, the fluid is thought of as a collection of particles whose path we trace as a function of time; the equations of motion roughly determine a map $\Phi_t(x)$ with $\Phi_0(x)=x$ determining the position at time $t\geq0$ of the particle located at $x$ when $t=0$.  Contrastingly, \emph{Eulerian} fluid dynamics takes the point of view of a barnacle attached to a point in the tank of water counting the number of particles that flow past a point $x$; this formulation might seek a function $\rho_t(x)$ giving the density of the fluid at a non-moving point $x$ as a function of time $t$.

So far, our formulation of transport has been Lagrangian:  The transportation plan $\pi$ explicitly determines how to match particles from the source distribution $\mu$ to the target distribution $\nu$.  Using a particularly clever change of variables, a landmark paper by Benamou \& Brenier shows that the 2-Wasserstein distance from~\eqref{eq:wasserstein} over Euclidean space with cost $c(x,y)=\|x-y\|_2^2$ can be computed in an Eulerian fashion~\cite{benamou2000computational}:
\begin{equation}\label{eq:bbnonconvex}
\W_2^2(\rho_0,\rho_1)\!=\!\left\{\!
\begin{array}{r@{\ }l} \min_{v(x,t),\rho(x,t)} & \frac{1}{2}\int_0^1\int_{\R^n} \rho(x,t)\,\|v(x,t)\|_2^2\,dA(x)\,dt\\
\st & \rho(x,0)\equiv \rho_0(x)\ \forall x\in\R^n\\
& \rho(x,1)\equiv \rho_1(x)\ \forall x\in\R^n\\
& \frac{\partial \rho(x,t)}{\partial t}=-\nabla\cdot (\rho(x,t)v(x,t))\\&\hspace{.2in}\forall x\in\R^n, t\in(0,1)
\end{array}
\right.
\end{equation}
Here, we assume that we are computing the 2-Wasserstein distance between two distribution functions $\rho_0(x)$ and $\rho_1(x)$.  This is often referred to as a \emph{dynamical} model of transport and can be extended to spaces like Riemannian manifolds~\cite{mccann2001polar}.

Formulation~\eqref{eq:bbnonconvex} comes with an intuitive physical interpretation.  The time-varying function $\rho(x,t)$ gives the density of a gas as a function of time $t\in[0,1]$, which starts out in configuration $\rho_0$ and ends in configuration $\rho_1$.  The constraint $\frac{\partial\rho}{\partial t}=-\nabla\cdot (\rho v)$ is the \emph{continuity equation}, which states that the vector field $v(x,t)$ is the \emph{velocity} of $\rho$ as it moves as a function of time while preserving mass.  Over all possible ways to ``animate'' the motion from $\rho_0$ to $\rho_1$, the objective function minimizes $\frac{1}{2}\rho\|v\|_2^2$ (mass times velocity squared):  the total kinetic energy!  

From a computational perspective, it can be convenient to replace velocity $v$ with momentum $J:=\rho \cdot v$ to obtain an equivalent formulation to~\eqref{eq:bbnonconvex}:
\begin{equation}\label{eq:bb}
\W_2^2(\rho_0,\rho_1)\!=\!\left\{\!
\begin{array}{r@{\ }l} \min_{J(x,t),\rho(x,t)} & \frac{1}{2}\int_0^1\int_{\R^n} \frac{\|J(x,t)\|_2^2}{\rho(x,t)}\,dA(x)\,dt\\
\st & \rho(x,0)\equiv \rho_0(x)\ \forall x\in\R^n\\
& \rho(x,1)\equiv \rho_1(x)\ \forall x\in\R^n\\
& \frac{\partial \rho(x,t)}{\partial t}=-\nabla\cdot J(x,t)\\&\hspace{.2in}\forall x\in\R^n, t\in(0,1)
\end{array}
\right.
\end{equation}
This formulation is convex jointly in the unknowns $(\rho,J)$.

\begin{figure}\centering
{\small\def\svgwidth{.5\textwidth}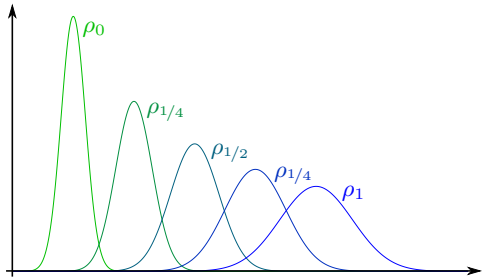}
\caption{Displacement interpolation from {\color[rgb]{0,.75,0}$\rho_0$} to {\color{blue}$\rho_1$} explains optimal transport between these two densities using a time-varying density function $\rho_t$, $t\in[0,1]$.}\label{fig:displacementinterpolation}
\end{figure}

Dynamical formulations of transport make explicit the phenomenon of \emph{displacement interpolation} \cite{mccann1994convexity,mccann1997convexity}, illustrated in Figure~\ref{fig:displacementinterpolation}.  Intuitively, the Wasserstein distance $\W_2$ between two distribution functions $\rho_0$ and $\rho_1$ is ``explained'' by a time-varying sequence of distributions $\rho_t$ interpolating from one to the other.  Unlike the trivial interpolation $\rho(t):=(1-t)\rho_0(x)+t\rho_1(x)$, optimal transport \emph{slides} the distribution across the geometric domain similar to a geodesic shortest path between points on a curved manifold.  Indeed, the intuitive connection to differential geometry is more than superficial:  \cite{otto2001geometry,lott2008some} show how to interpret~\eqref{eq:bbnonconvex} as a geodesic in an infinite-dimensional manifold of probability distributions over a fixed domain.

% mention w1/beckmann
Other $p$-Wasserstein distances $\W_p$ also admit Eulerian formulations.  Most importantly, the $1$-Wasserstein distance can be computed as follows:
\begin{equation}\label{eq:beckmann}
\W_1(\rho_0,\rho_1)=\left\{
\begin{array}{rl}
\min_{J(x)} & \int_{\R^n} \|J(x)\|_2\,dA(x)\\
\st & \nabla\cdot J(x)=\rho_1(x)-\rho_0(x).
\end{array}
\right.
\end{equation}
This problem, known as the \emph{Beckmann problem}, has connections to traffic modeling and other tasks in geometry.  From a computational perspective, it has the useful property that the vector field $J(x)$ has no time dependence, reducing the number of unknown variables in the optimization problem.

\section{Motivating Applications}\label{sec:app}

Having developed the basic definition and theoretical properties of the optimal transport problem, we can now divert from theoretical discussion to mention some concrete applications of transport in the computational world.  These are just a few, chosen for their diversity (and no doubt biased toward areas adjacent to the author's research); in reality optimal transport is beginning to appear in a huge variety of computational pipelines.  Our goal in this section is not to give the details of each problem and its resolution with transport, but just to give a flavor of how optimal transport can be applied as a powerful modeling tool in application-oriented disciplines as well as citations to more detailed treatments of each application.

\subsubsection*{Operations and logistics.}  Given its history and even its name, it comes as no surprise that a primary application of optimal transport is in the operations and logistics world, in which engineers are asked to find a minimum-cost routing of packages or materials to customers.  The basic theory and algorithms for this case of optimal transport date back to World War II, in which optimal transport of soldiers, weapons, supplies, and the like were by no means theoretical problems.  

A particular case of interest in this community is that of transport over a graph $G=(V,E)$.  Here, shortest-path distances over the edges of $G$ provide the costs for transport, leading to a problem known to computer scientists as \emph{minimum-cost flow without edge capacities}~\cite{ravindra1993network}.  This linear program is a classic algorithmic problem, with well-known algorithms including cycle canceling~\cite{klein1967primal}, network simplex~\cite{orlin1997polynomial}, and the Ford--Fulkerson method~\cite{ford1956solving}.  A challenge for theoretical computer scientists is to design algorithms achieving lower-bound time complexity for solving this problem; recent progress includes~\cite{sherman2017generalized}, which achieves a near-linear runtime using an approach that almost resembles a numerical algorithm rather than a discrete method.

\subsubsection*{Histogram-based descriptors.} Some of the earliest applications of optimal transport in the computational world come from computer vision~\cite{rubner2000earth}.  Suppose we wish to perform \emph{similarity search} on a database of photographs:  Given one photograph, we wish to search the database for other photos that look similar.  One reasonable way to do this is to describe each photograph as a histogram---or probability distribution---over the space of colors.  Two photographs roughly look similar if they have similar color histograms as measured using optimal transport distances (known in this community as the ``Earth Mover's Distance''), giving a simple technique for sorting and searching the dataset.

This basic approach comes up time and time again in the applied world.  For images, rather than binning colors into a histogram one could bin the orientations and strengths of the gradients to capture the distribution of edge features~\cite{pele2009fast}.  Recent work has proposed an embedding of the words in an English dictionary into Euclidean space $\R^n$~\cite{mikolov2013efficient}, in which case the words present in a given document become a point cloud or superposition of $\delta$-functions in $\R^n$; application of the Wasserstein (``Word-Movers'') distance in this case is an effective technique for document retrieval~\cite{kusner2015word}.

\subsubsection*{Registration.} Suppose we wish to use a medical imaging device such as the MRI to track the progress of a neurodegenerative disease.  On a regular basis, we might ask the subject to return to the laboratory for a brain scan, each time measuring a signal over the volume of the MRI indicating the presence or absence of brain tissue.  These signals can vary drastically from visit to visit, not just due to the progress of the disease but also due to more mundane issues like movement of the patient in the measurement device or nonrigid deformation of the brain itself. 

Inspired by issues like those mentioned above, the task of computing a map from one scan to another is known as \emph{registration}, and optimal transport has been proposed time and time again as a tool for this task.  The basic idea of these tools is to use the transport map $\pi$ as a natural way to transfer information from one scan to another~\cite{haker2004optimal}.  One caveat is worth highlighting:  Optimal transport does not penalize splitting mass or making non-elastic deformations in the optimal map, so long as points of mass individually do not move too far.  A few recent methods attempt to cope with this final issue, e.g.\ by combining transport with an elastic deformation method more common in medical imaging~\cite{feydy2017optimal} or by defining modified versions of optimal transport that are invariant to certain species of deformation~\cite{cohen1999earth,memoli2011gromov,solomon2016entropic}.% tannenbaum elastic warping maps, FX transport plus nonrigid

\begin{figure}\centering
\includegraphics[width=.3\linewidth]{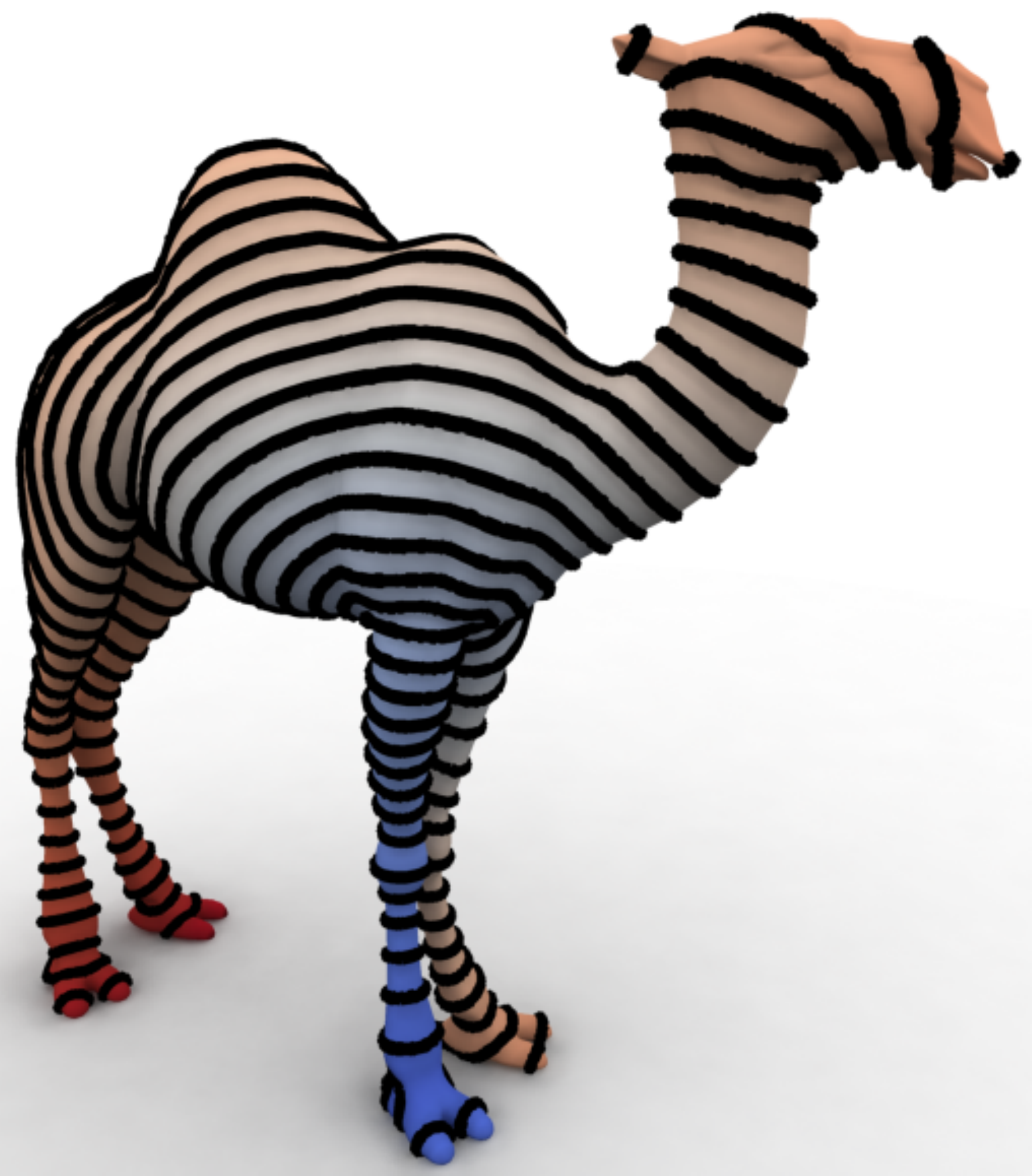}
\caption{Level sets of geodesic distance to the front right toe of a 3D camel model approximated using the optimal transport technique~\cite{solomon2014earth}.}\label{fig:levelset}
\end{figure}

\subsubsection*{Distance approximation.}  A predictable property of the $p$-Wasserstein distance $\W_p$ for distributions over a surface or manifold $\mathcal M$ is that the distance between $\delta$-functions centered at two points $x_0,x_1\in\mathcal M$ reproduces the geodesic (shortest-path) distance from $x_0$ to $x_1$.  While distances in Euclidean space are computable using a closed-form formula, distances along discretized surfaces can be challenging to compute algorithmically, requiring techniques like fast marching~\cite{sethian1999fast}, the theoretically-justified but difficult-to-implement MMP algorithm~\cite{mitchell1987discrete}, or diffusion-based approximation~\cite{crane2013geodesics}.  In this regime, fast algorithms for approximating optimal transport distances $\W_p$ restricted to $\delta$-functions actually provide a way to approximate geodesic distances while preserving the triangle inequality~\cite{solomon2014earth}; the level sets of one such approximation are shown in Figure~\ref{fig:levelset}.

\subsubsection*{Blue noise and stippling.} Certain laser printers and other devices can only print pages in black-and-white---no gray.  The idea of \emph{halftoning} is that gray values between black and white can be approximated in a perceptually reasonable fashion by patterns of black dots of varying radius or location over a white background; the halftoned image can be printed using the black-and-white printer and from a distance appears similar to the original.

\begin{figure}\centering
\includegraphics[width=\linewidth]{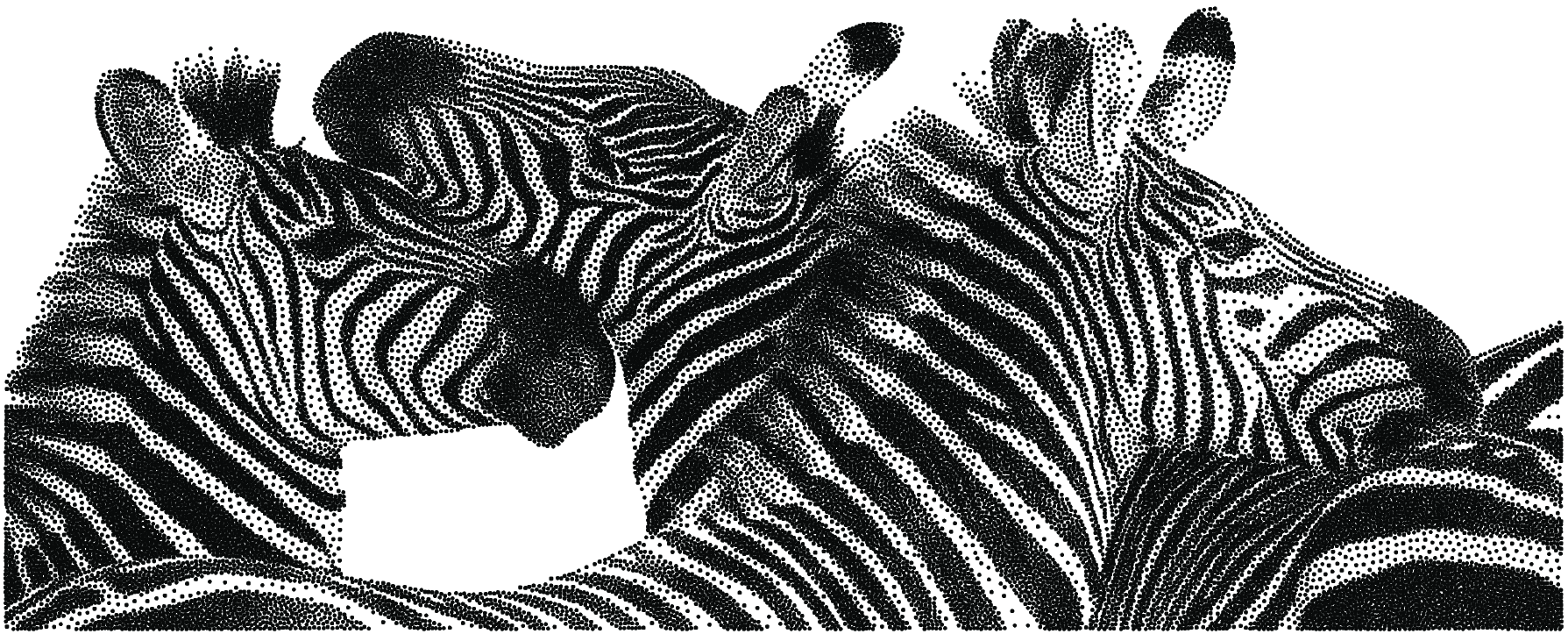}
\caption{A blue noise pattern generated using~\cite{de2012blue} (image courtesy F.\ de Goes, generated from photograph by F.\ Durand).}\label{fig:bluenoise}
\end{figure}
A reasonable model for halftoning involves optimal transport.  In particular, suppose we think of a grayscale image as a \emph{distribution} of ink on a white page; that is, the image can be understood as a measure $\mu\in\Prob([0,1]^2)$, where $[0,1]^2$ is the unit rectangle representing the image plane.  Under the reasonable assumption that ink is conserved, we might attempt to approximate $\mu$ with a set of dots of black inks, modeled using $\delta$-functions centered at $x_i$.  The intensity of the dot cannot be modulated (the printer only knows how to print in black-and-white), but the location can be moved, leading to an optimization problem to the effect:
$$
\min_{x_1,\ldots,x_n} \W_2^2\left( \mu, \frac{1}{n}\sum_i \delta_{x_i}\right).
$$
Here, the variables are the locations of the $n$ dots approximating the image, and the Wasserstein $2$-distance is used to measure how well the dots approximate $\mu$.  This basic idea is extended in~\cite{de2012blue} to a pipeline for computing \emph{blue noise}; an example of their output is shown in Figure~\ref{fig:bluenoise}.

\subsubsection*{Political redistricting.}  A few recent attempts to propose political redistricting procedures have incorporated ideas from optimal transport to varying degrees of success.  For example, optimal transport might provide one simplistic means of assigning voters to poling centers.   The distribution of voters over a map is ``transported'' to a sparse set of polling places, where distributional constraints reflect the fact that each polling center can only handle so many voters; assigning each voter to his/her closest polling center might cause polling centers in city centers to become overloaded. A few papers have proposed variations on this idea to design compact voting districts e.g.\ for the US House of Representatives~\cite{svec2007applying,miller2007problem,cohen2017balanced,Optimaldistricts}.  Many confounding---but incredibly important---factors obscure the application of this simplistic mathematical model in practice, ranging from compliance with civil rights law to the simple decision of a transport cost (e.g.\ geographical versus road network versus public transportation versus travel time).

\subsubsection*{Statistical estimation.} Parameter estimation is a key task in statistics that involves ``explaining'' a given dataset using a statistical model.  For example, given the set of heights of people in a room $\{h_1,\ldots,h_n\}$, a simple parameter estimation task might be to estimate the mean $h_0$ and standard deviation $\sigma$ of a normal (bell curve/Gaussian) distribution $g(h | h_0,\sigma)$ from which the data was likely sampled.  %minimum kantorovich estimator, bassetti 2006, MLE minimizes KL(data | p_theta); mention DRO and WASP

Principal among the techniques for parameter estimation is the \emph{maximum likelihood estimator} (MLE).  Continuing in our height data example, assuming the $n$ heights are drawn independently, the joint probability of observing the given set of heights in the room is given by the product
$$
P(h_1,\ldots,h_n|h_0,\sigma)=\prod_{i=1}^n g(h_i |h_0,\sigma).
$$
The MLE of the data is the estimate of $(h_0,\sigma)$ that maximizes this probability value:
$$
(h_0,\sigma)_{\mathrm{MLE}} := \arg\max_{h_0,\sigma} P(h_1,\ldots,h_n|h_0,\sigma).
$$
For algebraic reasons it is often easier to maximize the \emph{log likelihood} $\log P(\cdots),$ although this is obviously equivalent to the formulation above.

As an alternative to the MLE, however, the \emph{minimum Kantorovich estimator} (MKE)~\cite{bassetti2006minimum} uses machinery from optimal transport.  As the name suggests, the MKE estimates the parameters of a distribution by minimizing the transport distance between the parameterized distribution and the empirical distribution from data.  For our height problem, the optimization might look like
$$
(h_0,\sigma)_{\mathrm{MKE}} := \arg\min_{h_0,\sigma}\W_2^2\left(
\frac{1}{n}\sum_i \delta_{h_i}, g(\cdot|h_0,\sigma)
\right)
$$
The differences between MLE, MKE, and other alternatives can be subtle from the outside looking in, and the MKE is only recently being studied in applied environments in comparison to more conventional alternatives.  Since it takes into account the distance measure of the geometric space on which the samples are defined, the MKE appears to be robust to geometric noise that can confound more traditional alternatives---at the price of increased computational expense.  Recent applications have shown value of this estimator for training and inference in machine learning models~\cite{montavon2016wasserstein,bernton2017inference}.

\subsubsection*{Domain adaptation.}  Many basic statistical and machine learning algorithms make a false assumption that the ``training'' and ``test'' data are distributed equally.  As an example where this is not the case, suppose we wish to make an object recognition tool that learns how to label the contents of a photograph.  As training data, we use the listings on an e-commerce site like Amazon.com, which contain not only a photographs of a given object but also text describing it.  But, while this training data is extremely clean, it is not representative of possible test data, e.g.\ gathered by a robot navigating a shopping mall:  Photographs collected by the latter likely contain clutter, a variety of lighting configurations, and countless other confounding factors.  Algorithms designed to compensate for the difference between training and test data are known as \emph{domain adaptation} techniques.

One possibility is to use optimal transport to design a stable domain adaptation tool.  The basic idea is to view the training data as a point cloud in some Euclidean space $\R^d$.  For instance, perhaps $d$ could be the number of pixels in a photograph; the location of every point in the point cloud determines the contents of the photo, and as additional information each point is labeled with a text name.  The test data is also a point cloud in $\R^d$, but thanks to the confounding factors listed above perhaps these two points clouds are not aligned.  \cite{courty2017optimal} proposes using optimal transport to align the training data to the test data and to carry the label information along, e.g.\ attempting to align the space of Amazon.com photos to the space of shopping mall photos.  Once the training and test data are aligned, it makes sense to transfer information, classifiers, and the like from one to the other.

\begin{figure}\centering
\includegraphics[height=1.2in]{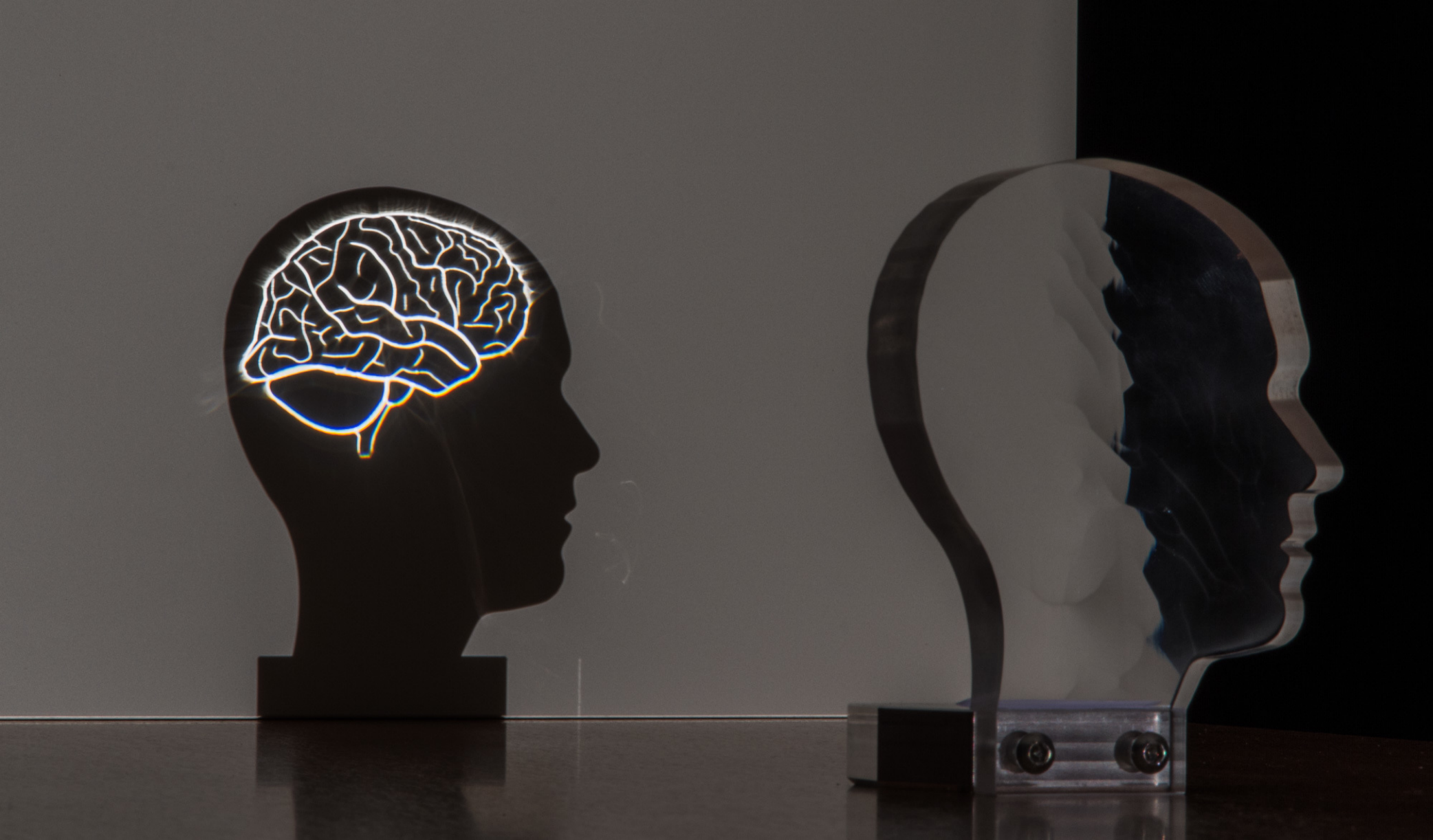}
\includegraphics[height=1.2in]{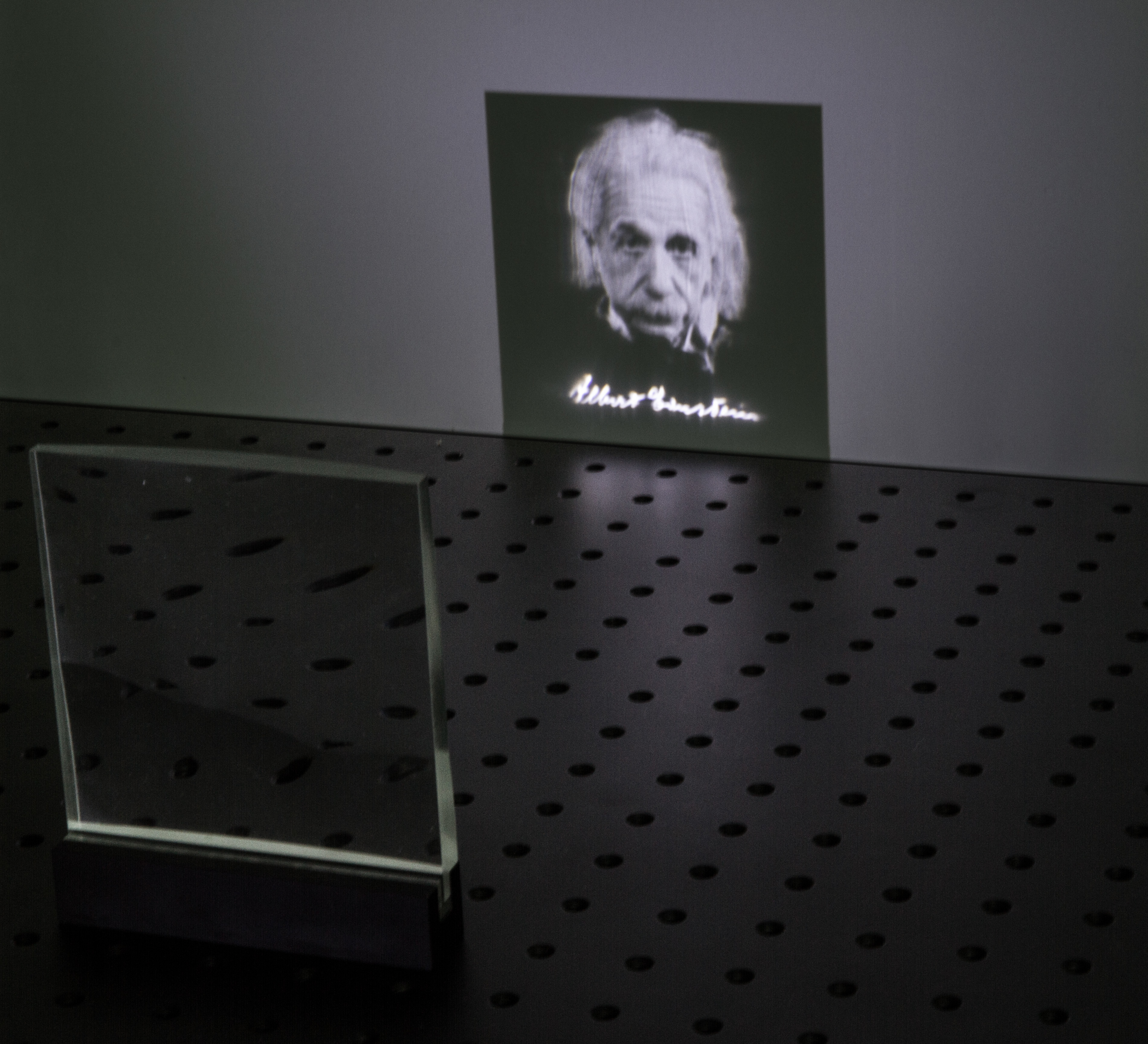}
\includegraphics[height=1.2in]{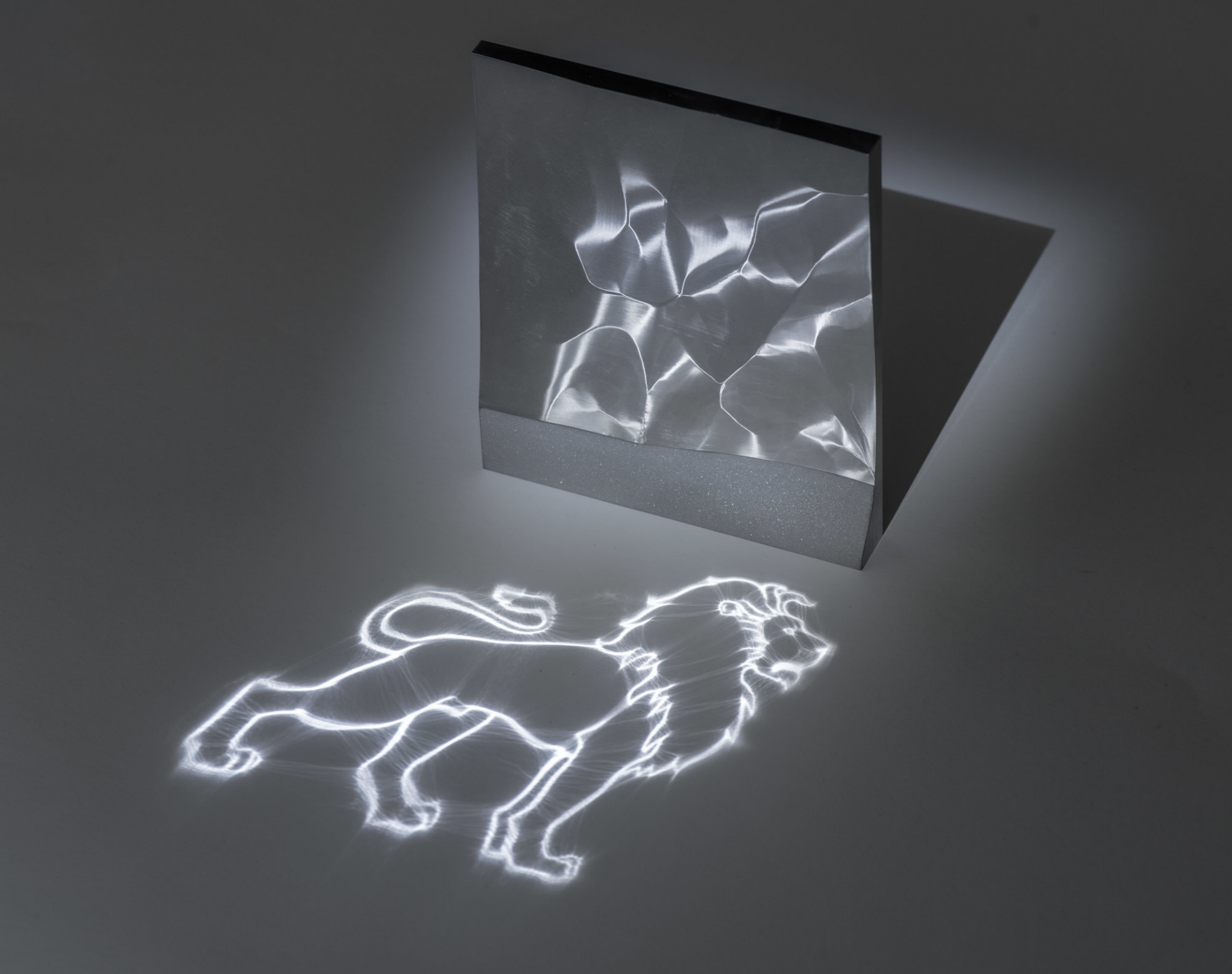}
\caption{Optimal transport is used to design the shape of transparent or reflective material to show a particular caustic pattern (image courtesy of EPFL Computer Graphics and Geometry Laboratory and Rayform SA).}\label{fig:caustic}
\end{figure}

\subsubsection*{Engineering design.}  Optimal transport has found application in design tools for many engineering tasks, from reflector design~\cite{oliker1987near,wang1996design} to aerodynamics~\cite{plakhov2012billiards}.  One intriguing paper uses optimal transport to design transparent objects made of materials like glass, which can focus light into \emph{caustics} via refraction~\cite{schwartzburg2014high}.  By minimizing the transport distance between the light rays by the glass and a desired black-and-white image, they can ``shape'' the distribution of light as it comes out of a window.  An example caustic design computed using their method is shown in Figure~\ref{fig:caustic}.  % n. Oliker [72] and X-J Wang [90] have pioneered the use of transportation theory in reflector design, while Plakhov has been exploring novel applications in aerodynamics, see e.g., [76, 77]. Also caustics

\section{One Problem, Many Discretizations}\label{sec:discrete}

%desiderata: structure preservation, efficiency, stability

Computational optimal transport is a relatively new discipline, and techniques for solving the optimal transport problem and in particular computing Wasserstein distances are still a topic of active research.  So far, it appears that no ``one size fits all'' approach has been discovered; rather, different applications and scenarios demand different numerical techniques for optimal transport.

Several desiderata inform the design of an algorithm for optimal transport:
\begin{itemize}[leftmargin=*]
\item \textsc{Efficiency:}  While $L_1$ distances and KL divergences are computable using closed-form formulas, optimal transport distance computation requires solving an optimization problem.  The cost of solving this problem relative to the cost of direct computation of transport's simpler alternatives is largely the reason why optimal transport has not reached a higher level of popularity in the applied world.  But this scenario is changing:  New high-speed algorithms for large-scale transport are nearly competitive with more traditional alternatives while bringing to the table the geometric structure unique to transport world.
\item \textsc{Stability:}  A theme in the numerical analysis literature is stability, the resilience of a computation to small changes in the input.  Stability of the minimal transport objective value and/or its accompanying transport map can be a challenging topic.  Linear program discretizations of continuum optimal transport problems tend to resemble~\eqref{eq:discrete_ot} above, a linear program whose optimal solution $T$ \emph{provably} has the sparsity of a permutation matrix; this implies that a small perturbation of $v$ or $w$ may result in a discrete change of $T$'s sparsity structure.
\item \textsc{Structure preservation:}  Transport is well-studied theoretically, and one could reasonably expect that key properties of transport in the infinite-dimensional case are preserved either exactly or approximately when they are computed numerically.  For instance, Wasserstein distances enjoy a triangle inequality, and Eulerian formulations of transport have connections to gradient flows and other PDE.  Provable guarantees that these structures are preserved in discretizations of transport help assure that nothing critical is lost in the process of approximating transport distances algorithmically.
\end{itemize}

One reason why there are so many varied algorithms available for numerical OT is that the problem can be written in so many different ways (see \S\ref{sec:manyformulas}).  A basic recipe for designing a transport algorithm is to choose any one of many equivalent formulations of transport---all of which yield the same optimal value in theory---, discretize any variables that are otherwise infinite-dimensional, and design a bespoke optimization algorithm to solve the resulting problem, which now has a finite number of variables.  The flexibility of choosing \emph{which} version of transport to discretize usually is tuned to the geometry of a given application, desired properties of the resulting discretization, and ease of optimizing the discretized problem.  The reality of choosing a discretization to facilitate ease of computation reflects a tried-and-true maxim of engineering:  ``If a problem is difficult to solve, change the problem.''

In this section, we roughly outline a few discretizations and accompanying optimization algorithms for numerical OT. Our goal is not to review all well-known techniques for computational transport thoroughly but rather to highlight the breadth of possible approaches and to give a few practical pointers for implementing state-of-the-art transport algorithms at home.

\subsection{Regularized Transport}

We will start by describing \emph{entropically-regularized transport}, a technique that has piqued the interest of the machine learning community after its introduction there in 2013~\cite{cuturi2013sinkhorn}.  This technique has an explicit trade-off between accuracy and computational efficiency and has shown particularly strong promise in the regime where a rough estimate of transport is sufficient.  This regime aligns well with the demands of ``big data'' applications, in which individual data points are likely to be noisy---so obtaining an extremely accurate transport value would be overkill computationally.

Regularization is a key technique in optimization and inverse problems in which an objective function is modified to encode additional assumptions and/or to make it easier to minimize.  For example, suppose we wish to solve the least-squares problem $\min_x \|Ax-b\|_2^2$ for some $A\in\R^{m\times n}$ and $b\in\R^m$.  When $A$ is rank-deficient or if $m<n$, an entire affine space of $x$'s achieve the minimal value.  To get around this, we could apply Tikhonov regularization (also known as ridge regression), in which we instead minimize $\|Ax-b\|_2^2+\alpha\|x\|_2^2$ for some $\alpha>0$.  As $\alpha\rightarrow0$ a solution of the original least-squares problem is recovered, while for any $\alpha>0$ the regularized problem is guaranteed to have a unique minimizer; as $\alpha\rightarrow\infty$, we have $x\rightarrow0$, a predictable but uninteresting value.  From a high level, we can think of $\alpha$ as trading off between fidelity to the original problem $Ax\approx b$ and ease of solution:  For small $\alpha>0$ the problem is near-singular but close to the original least-squares formulation, while larger $\alpha$ makes the problem easier to solve.

The variables in the basic formulation of transport are nonnegative probability values, which do not appear to be amenable to standard least-squares style Tikhonov regularization.  Instead, entropic regularization uses a regularizer from information theory:  the entropy of a probability distribution.  Suppose a probability measure has distribution function $\rho(x)$.  The (differential) entropy of $\rho$ is defined as
\begin{equation}\label{eq:entropy}
H[\rho]:=-\int \rho(x)\log\rho(x)\,dx.
\end{equation}
This definition makes two assumptions that are needed to work with entropy, that a probability measure admits a distribution and that it is nonzero everywhere---otherwise $\log\rho(x)$ is undefined.  $H[\rho]$ is a concave function of $\rho$ that roughly measures the ``fuzziness'' of a distribution.  Low entropy indicates that a distribution is sharply peaked about a few points, while high entropy indicates that it is more uniformly distributed throughout space.

The basic approach in entropically-regularized transport is to add a small multiple of $-H[\pi]$ to regularize the transport plan $\pi$ in the OT problem.  We will start by discussing the discrete problem~\eqref{eq:discrete_ot}, which after entropic regularization can be written as follows:
\begin{equation}\label{eq:entropic_ot}
\OT_\alpha(v,w;C):=\left\{
\begin{array}{rl}
\min_{T\in\R^{k_1\times k_2}} & \sum_{ij} T_{ij}c_{ij}+\alpha\sum_{ij} T_{ij}\log T_{ij}\\
\st 
& \sum_j T_{ij}=v_i\ \forall i\in\{1,\ldots,k_1\}\\
& \sum_i T_{ij}=w_j\ \forall j\in\{1,\ldots,k_2\}.
\end{array}
\right.
\end{equation}
We are able to drop the $T\geq0$ constraint because $\log T_{ij}$ in the objective function prevents negative $T$ values.  

The objective function from~\eqref{eq:entropic_ot} can be refactored slightly:
\begin{align}
\sum_{ij} T_{ij}c_{ij}+\alpha\sum_{ij} T_{ij}\log T_{ij}
&=\alpha\sum_{ij} T_{ij} \left(\frac{c_{ij}}{\alpha}+\log T_{ij}\right)\nonumber\\
&=\alpha \sum_{ij} T_{ij} \log \frac{T_{ij}}{e^{-c_{ij}/\alpha}}\nonumber\\
&=\alpha \KL(T|K_\alpha).\label{eq:kl}
\end{align}
Here, we define a \emph{kernel} $K_\alpha$ via $(K_\alpha)_{ij}:=e^{-c_{ij}/\alpha}$. $\KL$ denotes the Kullback--Leibler divergence~\cite{kullback1951information}, a distance-like (but asymmetric) measure of the similarity between $T$ and $K$ from information theory; the definition of $K_\alpha$ is singular when $\alpha=0$, indicating that the connection to $\KL$ is only possible in the $\alpha>0$ regime.  

\begin{figure}\centering
\begin{tabular}{c@{\hspace{.75in}}c}
{\small\def\svgwidth{.4\textwidth}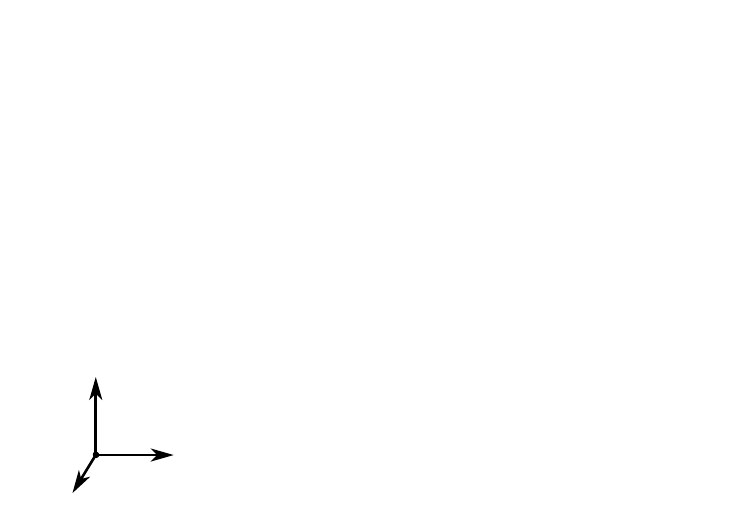}&
{\small\def\svgwidth{.4\textwidth}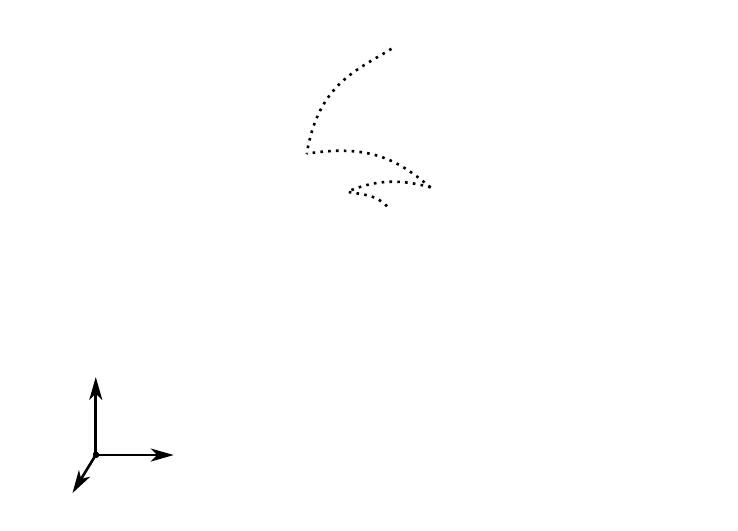}\\
(a) Projection & (b) Alternating projection
\end{tabular}
\caption{(a) Intuition for the optimization problem~\eqref{eq:entropic_ot} as a projection of $K_\alpha$ onto the prescribed {\color{red} row sum} and {\color{blue} column sum} constraints with respect to KL divergence~\eqref{eq:kl}. (b) The Sinkhorn algorithm projects back and forth onto one set of constraints and then the other, converging to the transport matrix $T^\ast$.}\label{fig:klprojection}
\end{figure}

Equation~\eqref{eq:kl} gives an intuitive explanation for entropy-regularized transport illustrated in Figure~\ref{fig:klprojection}(a).  The matrix $K$ does not satisfy the constraints of the regularized transport problem~\eqref{eq:entropic_ot}.  Thinking of $\KL$ roughly as a distance measure, our job is to find the \emph{closest projection} (with respect to $\KL$) of $K$ onto the set of $T$'s satisfying the constraints $\sum_j T_{ij}=v_i$ and $\sum_i T_{ij}=w_j$.  With this picture in mind, Figure~\ref{fig:klprojection}(b) illustrates the Sinkhorn algorithm for entropy-regularized transport derived below, which alternates between projecting onto one of these sets and then the other.

Continuing in our derivation, we return to~\eqref{eq:entropic_ot} to derive first-order optimality conditions.  Since~\eqref{eq:entropic_ot} is an equality-constrained differentiable minimization problem, it can be solved using a standard multi-variable calculus technique:  the method of Lagrange multipliers.  There are $k_1+k_2$ constraints, so we need $k_1+k_2$ Lagrange multipliers, which---following the derivation of~\eqref{eq:discrete_dual}---we store in vectors $\phi\in\R^{k_1}$ and $\psi\in\R^{k_2}$.  The Lagrange multiplier function here is:
\begin{align*}
\Lambda(T;\phi,\psi)
&:=
\sum_{ij} T_{ij}c_{ij}+\alpha\sum_{ij} T_{ij}\log T_{ij}
\\&\hspace{.25in}+
\sum_i \phi_i \left(v_i-\sum_j T_{ij}\right)
+
\sum_j \psi_j \left(w_j-\sum_i T_{ij}\right)
\\
&=\langle T,C\rangle + \alpha \langle T,\log T\rangle + \phi^\top (v-T\1) + \psi^\top(w-T^\top \1)
\end{align*}
Here, $\langle\cdot,\cdot\rangle$ indicates the element-wise inner product of matrices, the log is element-wise, and $\1$ indicates the vector of all ones.  Taking the gradient with respect to $T$ gives the following first-order optimality condition:\footnote{Readers uncomfortable with this sort of calculation are strongly encouraged to take a look at the useful ``cheat sheet'' document~\cite{petersen2008matrix}.}
\begin{align*}
&0=\nabla_T\Lambda
= C + \alpha \1\1^\top + \alpha\log T - \phi\1^\top - \1\psi^\top \\
&\implies \log T 
= \frac{(\phi-\alpha\1)\1^\top}{\alpha} + \frac{\1\psi^\top}{\alpha} +\log K_\alpha\textrm{ where }K_\alpha :=\exp[-C/\alpha]\\
&\implies \boxed{T
=\diag[p] K_\alpha \diag[q]}\textrm{ where }p:=\exp\left[\frac{\phi-\alpha\1}{\alpha}\right]\textrm{ and }
q:=\exp\left[\frac{\psi}{\alpha}\right].
\end{align*}
Here, $\diag[v]$ indicates the diagonal matrix whose diagonal is $v$.  The key result is the boxed equation, which gives a formula for the unknown transport matrix $T$ in terms of two unknown vectors $p$ and $q$ derived by changing variables from the Lagrange multipliers $\phi$ and $\psi$.  There are multiple choices of $p$ and $q$ in terms of $\phi$ and $\psi$ that all give the same ``diagonal rescaling'' formula including some that are more symmetric, but this detail is not important.

Next we plug the new relationship $T=\diag[p] K_\alpha \diag[q]$ into the constraints of~\eqref{eq:entropic_ot} to find
\begin{equation}\label{eq:sinkhorn_relationships}
\begin{array}{r@{\ }l}
p\otimes (K_\alpha q)&=v\\
q\otimes (K_\alpha^\top p)&=w.
\end{array}
\end{equation}
Here, $\otimes$ denotes the elementwise (Hadamard) product of two vectors or matrices.  These formulas determine the unknown vector $p$ in terms of $q$ and vice versa.

The formulas~\eqref{eq:sinkhorn_relationships} directly suggest a state-of-the-art technique for entropy-regularized optimal transport, known as the \emph{Sinkhorn (or Sinkhorn--Knopp) algorithm} and dating back to an early technique for matrix rescaling~\cite{sinkhorn1967concerning}.  This extremely succinct algorithm successively updates estimates of $p$ and $q$. Iteration $k$ is given by the update formulas ($\oslash$ denotes elementwise division)
\begin{align*}
p^{k+1}&\gets v\oslash (K_\alpha q^k)\\
q^{k+1}&\gets w\oslash (K_\alpha^\top p^{k+1}).
\end{align*}
It can be implemented in fewer than ten lines of code!  The basic approach is to update $p$ in terms of $q$ using the first relationship, then $q$ in terms of $p$ using the second relationship, then $p$ again, and so on.  Using essentially the geometric intuition provided in Figure~\ref{fig:klprojection}(b) for this technique and explored in-depth in~\cite{benamou2015iterative}, one can prove that $\diag[p] K_\alpha \diag[q]$ converges asymptotically to the optimal $T$ at a relatively efficient rate regardless of the initial guess.

% mention parallel computation, fast linear algebra
Several advantages distinguish the Sinkhorn method from its peers in the numerical optimization world.  Most critically, beyond its ease of implementation, this algorithm is built from simple linear algebra operations---matrix-vector multiplies and elementwise arithmetic---that parallelize well and can be carried out extremely quickly on modern processing hardware.  One modern spin on Sinkhorn shows how to shave off even more calculations while preserving its favorable convergence rate~\cite{altschuler2017}.

Beyond inspiring a huge body of follow-on work in machine learning and computer vision, the Sinkhorn rescaling algorithm provides a means to adapt optimal transport to discrete domains suggested in~\cite{solomon2015convolutional}.  So far, our description of the Sinkhorn method has been generic to \emph{any} cost matrix $C$. % and (2) discrete rather than applying to transport over smooth manifolds.  
Adding geometric structure to the problem gives it a strong interpretation using heat flow and suggests a faster way to carry out Sinkhorn iterations on discrete domains.  %Since our survey involves discrete domains anyway, we will not address (2) directly but will skip directly to a discrete formulation on surfaces like triangle meshes that also suggests how to make Sinkhorn work in the continuum case for theoretical purposes.

%We will study the case illustrated in Figure \addreference.  
Suppose that the transport cost $C$ is given by squared pairwise distances along a discretized piece of geometry such as a triangulated surface, denoted $\Sigma$; this corresponds to computing a regularized version of the 2-Wasserstein distance~\eqref{eq:wasserstein}.  The dual variables $p$ and $q$ can be thought of as \emph{functions} over $\Sigma$, discretized e.g.\ using one value per vertex.  Then, the kernel $K_\alpha$ has elements $$(K_\alpha)_{ij}=e^{-d(x_i,x_j)^2/\alpha},$$ where $d(x_i,x_j)$ denotes the shortest-path (geodesic) distance along the domain from vertex $i$ to vertex $j$.

To start, if our domain is flat, or Euclidean, then $(K_\alpha)_{ij}=e^{-\|x_i-x_j\|_2^2/\alpha}$ for points $\{x_i\}_i\subseteq \R^n$.  Considered as a function of the $x_i$'s, we recognize $K_\alpha$ up to scale as a \emph{Gaussian} (or normal distribution, or bell curve) in distance.  Multiplication by $K_\alpha$ is then \emph{Gaussian convolution}, an extremely simple operation that can be carried out algorithmically using methods like the Fast Fourier Transform (FFT).  In other words, rather than explicitly computing and storing the matrix $K_\alpha$ as an initial step and computing matrix-vector products $K_\alpha p$ and $K_\alpha q$ (note $K_\alpha$ is symmetric in this case) in every iteration of the Sinkhorn algorithm, in this case we can replace these products with convolutions $g_\sigma \ast p$ and $g_\sigma \ast q$, where $\ast$ denotes convolution and $g_\sigma$ is a Gaussian whose standard deviation is determined by the regularizer $\alpha$.  This is \emph{completely equivalent} to the Sinkhorn method that explicitly computes the matrix-vector product, while eliminating the need to store $K_\alpha$ and improving algorithmic speed thanks to fast Gaussian convolution.  Put more simply, in the Euclidean case \textbf{multiplication by $K_\alpha$ is more efficient than storing $K_\alpha$} since we can carry out the former implicitly.

When $\Sigma$ is curved, we can use a mathematical sleight of hand modifying the entropic regularizer to improve computational properties while maintaining convergence to the true optimal transport value as the regularizer goes to zero.  We employ a well-known property of geodesic distances introduced in theory in~\cite{varadhan1967behavior} and applied to computing distances on discrete domains in~\cite{crane2013geodesics}.  This property, known as Varadhan's formula, states that geodesic distance $d(x,y)$ between two points $x,y$ on a manifold can be recovered from heat diffusion over a short time:
$$
d(x,y)^2=\lim_{t\rightarrow0}[-2t\ln\mathcal H_t(x,y)].
$$
Recall that the heat kernel $\mathcal H_t(x,y)$ determines diffusion between $x,y\in \Sigma$ after time $t$.  That is, if $f_t$ satisfies the heat equation $\partial_tf_t=\Delta f_t$, where $\Delta$ denotes the Laplacian operator, then
$$
f_t(x)=\int_\Sigma f_0(y)\mathcal H_t(x,y)\,dy.
$$

Connecting to the previous paragraph, the heat kernel in Euclidean space is exactly the Gaussian function!  Hence, if we replace the kernel $K_\alpha$ with the heat kernel $\mathcal H_{\alpha/2}$ in Sinkhorn's method, in the Euclidean case nothing has changed.  In the curved case, we get a new approximation of Wasserstein distances introduced as ``convolutional Wasserstein distances.''

All that remains is to convince ourselves that we can compute matrix-vector products $\mathcal H_t \cdot p$ when $\mathcal H_t$ is the heat kernel of a discretized domain $\Sigma$ that is not Euclidean.  Thankfully, armed with material from other chapters in this tutorial, this is quite straightforward in the context of discrete differential geometry.  In particular, the well-known cotangent approximation of the Laplacian $\Delta$ can be combined with standard ordinary differential equation (ODE) solution techniques to carry out heat diffusion in this case using sparse linear algebra.  We refer the reader for~\cite{solomon2015convolutional} for details of one implementation that uses DDG tools extensively.

\subsection{Eulerian Algorithms}

\begin{figure}
\centering
{\small\def\svgwidth{.9\textwidth}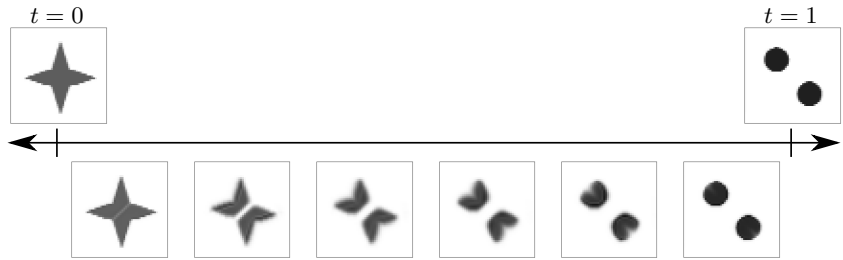}
\caption{Output from an Eulerian algorithm for optimal transport extending~\cite{benamou2000computational} (image courtesy H.\ Lavenant); interpolation between the two distributions on the top is shown below the timeline.  In addition to finding the transport cost, methods in this class also provide a sequence of distributions interpolating between the two inputs.}\label{fig:bb_example}
\end{figure}

Entropically-regularized transport works with the Kantorovich formulation~\eqref{eq:transport}.  This may be one of the earliest and most intuitive definitions of optimal transport, but this in itself is not a strong argument in favor of tackling this formulation numerically.  As a point of contrast, we now explore a \emph{completely different} approximation of Wasserstein distances that can be useful in low-dimensional settings, built from the Eulerian (fluid mechanics) formulation of the 2-Wasserstein distance $\W_2^2$~\eqref{eq:bbnonconvex}.  Historically, this method pre-dates the popularity of entropically-regularized transport and has distinct advanges and disadvantages:  It explicitly computes a time-varying displacement interpolation of a density ``explaining'' the transport (see Figure~\ref{fig:bb_example}) but in the process must solve a difficult boundary-value PDE problem.  Beyond the original paper~\cite{benamou2000computational}, we recommend the excellent tutorial~\cite{peyre2010optimal} that steps through an implementation of this technique in practice.%; we use their modified version of~\cite{benamou2000computational} in our discussion below.

%Once again, when our problem becomes hard to discretize, we change the problem!  In particular, we make one 
We make a few more simplifications to the continuum formulation before discretizing it. 
%Just like our derivation of the dual formula in \S\ref{sec:manyformulas}, we can derive first-order optimality conditions for~\eqref{eq:bb} by writing it as a minimax problem
%$$
%\inf_{J,\rho}\sup_\phi \left\{ \int_0^1 \int_{\R^n} \left[
%\frac{\|J(x,t)\|^2}{\rho(x,t)} - \phi(x,t)\left(
%\frac{\partial \rho(x,t)}{\partial t}+\nabla\cdot J(x,t)
%\right)
%\right]\,dA(x)\,dt
%-
%\int_{\R^n} [\phi(x,0)\rho_0(x)-\phi(x,1)\rho_1(x)]\,dA(x)\right\}
%$$
%Here, $\phi(x,t)$ is a dual variable used to enforce the constraint $\frac{\partial\rho}{\partial t}=-\nabla\cdot J$ for $t\in(0,1)$ as well as the boundary conditions $\rho(x,0)=\rho_0(x)$ and $\rho(x,1)=\rho_1(x)$ at $t=0$ and $t=1$, respectively.
We start by making a quick observation:  for any vector $J\in\R^n$ and $\rho>0$ we have
$$
\frac{\|J\|_2^2}{2\rho}=\left\{\begin{array}{rl}
\sup_{a\in\R,b\in\R^n} & a\rho+b^\top J\\
\st & a+\frac{\|b\|_2^2}{2}\leq0.
\end{array}\right.
$$
This convex program not only justifies that the quotient $\nicefrac{\|J\|_2^2}{2\rho}$ is convex jointly in $J$ and $\rho$, but also it shows we can write the optimization problem~\eqref{eq:bb} with additional variables as
$$
\begin{array}{rl} \inf_{J,\rho}\sup_{a, b} & \int_0^1\int_{\R^n} [a(x,t)\rho(x,t)+b(x,t)^\top J(x,t)]\,dA(x)\,dt\\
\st & \rho(x,0)\equiv \rho_0(x)\ \forall x\in\R^n\\
& \rho(x,1)\equiv \rho_1(x)\ \forall x\in\R^n\\
& \frac{\partial \rho(x,t)}{\partial t}=-\nabla\cdot J(x,t)\ \forall x\in\R^n, t\in(0,1)\\
& a(x,t)+\frac{\|b(x,t)\|_2^2}{2}\leq0\ \forall x\in\R^n, t\in(0,1).
\end{array}
$$
Next, we introduce a dual potential function $\phi(x,t)$ similarly to the derivation of~\eqref{eq:transportdual} to take care of all but the last constraint:
\begin{equation}\label{eq:bb_partial}
\begin{array}{r@{\ }l} \inf_{J,\rho}\sup_{a, b,\phi} & \int_0^1\!\int_{\R^n}\!\Big[a(x,t)\rho(x,t)+b(x,t)^\top J(x,t)\!\\&\hspace{.75in}+\phi(x,t)\left(\frac{\partial\rho(x,t)}{\partial t}\!+\!\nabla\cdot J(x,t)\right)\Big]\,dA(x)\,dt\\
&+
\int_{\R^n} [\phi(x,\!1)(\rho_1(x)\!-\!\rho(x,\!1))\!-\!\phi(x,0)(\rho_0(x)\!-\!\rho(x,\!0))]\,dA(x)\\
\st & a(x,t)+\frac{\|b(x,t)\|_2^2}{2}\leq0\ \forall x\in\R^n, t\in(0,1).
\end{array}
\end{equation}
We can simplify some terms in this expression.  First, using integration by parts we have
$$
\int_0^1 \phi(x,t)\frac{\partial\rho(x,t)}{\partial t}\,dt
=
[\rho(x,1)\phi(x,1)-\rho(x,0)\phi(x,0)]-\int_0^1 \rho(x,t)\frac{\partial\phi(x,t)}{\partial t}\,dt
$$
%This  yields a new objective function
%$$
%\int_{\R^n} 
%\left\{
%\int_0^1
%\left(\rho(x,t)\left[a(x,t) - \frac{\partial\phi(x,t)}{\partial t}\right]+b(x,t)^\top J(x,t)+\phi(x,t)\nabla\cdot J(x,t)\right)
%\,dt
%-\phi(x,0)\rho_0(x)+\phi(x,1)\rho_1(x))
%\right\}\,dA(x)
%$$
We also can integrate by parts in $x$ to show
$$
\int_{\R^n} \phi(x,t)\nabla\cdot J(x,t)\,dA(x)=-\int_{\R^n} J(x,t)^\top \nabla\phi(x,t)\,dA(x).
$$
This simplification works equally well if we replace $\R^n$ with the box $[0,1]^n$ with periodic boundary conditions.  Incorporating these two integration by parts formulae into our objective function yields a new one:
\begin{align*}
\int_{\R^n} 
\Bigg\{
\int_0^1&
\left(\rho(x,t)\left[a(x,t) - \frac{\partial\phi(x,t)}{\partial t}\right]+J(x,t)^\top [b(x,t)-\nabla \phi(x,t)]\right)
\,dt\\&
-\phi(x,0)\rho_0(x)+\phi(x,1)\rho_1(x))
\Bigg\}\,dA(x)
\end{align*}
We now make some notational simplifications.  Define $z:=\{\rho,J\}$ and $q:=\{a,b\}$ with inner product $$\langle z,q\rangle:=\int_{\R^n}\int_0^1 (a(x,t)\rho(x,t)+b(x,t)^\top J(x,t))\,dt\,dA(x).$$
Furthermore, define
\begin{align*}
F(q)&:=\left\{
\begin{array}{ll}
0 & \textrm{ if } a(x,t)+\frac{\|b(x,t)\|_2^2}{2}\leq0\ \forall x\in\R^n, t\in(0,1)\\
\infty & \textrm{ otherwise.}
\end{array}
\right.\\
G(\phi) &:=
\int_{\R^n}
(\phi(x,0)\rho_0(x)-\phi(x,1)\rho_1(x))
\,dA(x)
\end{align*}
These functions are both convex.  These functions, plus our simplifications and a sign change, allow us to write~\eqref{eq:bb_partial} in a compact fashion as:
\begin{equation}\label{eq:bb_slick}
-\sup_z\inf_{q,\phi} \left[
F(q) + G(\phi) + \langle z, \nabla_{x,t}\phi -q  \rangle
\right],
\end{equation}
where $\nabla_{x,t}\phi:=\{\nicefrac{\partial\phi}{\partial t},\nabla_x\phi\}.$  

Blithely assuming strong duality, namely that we can swap the supremum and the infimum, we arrive at an alternative interpretation of~\eqref{eq:bb_slick}.  In particular, we can view $z$ as a Lagrange multiplier corresponding to a constraint $q=\nabla_{x,t}\phi$.  From this perspective, we actually can find a saddle point (max in $z$, minimum in $(q,\phi)$) of the \emph{augmented Lagrangian} $L_r$ for any $r\geq 0$:
$$
L_r(\phi,q,z):=F(q)+G(\phi)+\langle z, \nabla_{x,t}\phi -q \rangle+\frac{r}{2}\langle \nabla_{x,t}\phi -q, \nabla_{x,t}\phi  -q \rangle.
$$
The extra term here effectively adds zero to the objective function, assuming the constraint is satisfied.

The algorithm proposed in~\cite{benamou2000computational} iteratively updates estimates $(\phi^\ell,q^\ell,z^\ell)$ by cycling through the following three steps:
\begin{align*}
\phi^{\ell+1} &\gets \arg\min_\phi L_r(\phi,q^\ell,z^\ell)\\
q^{\ell+1} &\gets \arg\min_q L_r(\phi^{\ell+1},q,z^\ell)\\
z^{\ell+1} &\gets z^\ell + r(q^{\ell+1}-\nabla_{x,t}\phi^{\ell+1}).
\end{align*}
The first two steps update some variables while holding the rest fixed to the best possible value.  The third step is gradient step for $z$.  This cycling algorithm and equivalent formulations has many names in the literature---including ADMM~\cite{boyd2011distributed}, the Douglas--Rachford algorithm~\cite{douglas1956numerical,lions1979splitting}, and the Uzawa algorithm~\cite{uzawa196810}---and is known to converge under weak assumptions.

The advantage of this algorithm is that the individual update formulae are straightforward.  In particular, the $\phi$ update is equivalent to solving a Laplace equation $$\Delta_{x,t} \phi^{\ell+1}=\nabla_{x,t}\cdot(z^\ell-rq^\ell),$$
where $\Delta_{x,t}$ is the Laplacian operator in time and space.  The $q$ update decouples over $x$ and $t$, amounting to projecting $\nabla_{x,t}\phi^{\ell+1}+\nicefrac{z^\ell}{r}$ onto the constraints in the definition of $F(q)$ with respect to $L_2$, a one-dimensional problem solvable analytically.  And, the $z$ update is already in closed-form.

So far, we have described the Benamou--Brenier algorithm using continuum variables, but of course at the end of the day we must discretize the problem for computational purposes.  The most straightforward discretization assumes $\rho_0$ and $\rho_1$ are supported in the unit square $[0,1]^n$, which is broken up into a $m\times m\times\cdots\times m$ grid, and further discretizes the time variable $t\in[0,1]$ into $p$ steps.  Then, all degrees of freedom $(\phi,q,z)$ can be put on the grid vertices and interpolated in between using multilinear basis functions; this leads to a finite element (FEM) discretization of the problem that can be approached using techniques discussed in earlier chapters.    An alternative grid-based discretization and accompanying optimization algorithm is also given in~\cite{papadakis2014optimal}.

The use of PDE language makes this dynamical formulation of transport seem attractive as potentially compatible with machinery like discrete exterior calculus (DEC)~\cite{hirani2003discrete}, which could define a discrete notion of transport on simplicial complexes like triangle meshes that discretize curved surfaces.  This remains an open problem for challenging technical reasons.\footnote{Interested readers are encouraged to contact the author of this tutorial for preliminary results on this problem!}  Principally, discretizing the objective function $\nicefrac{\|J\|_2^2}{\rho}$ on a triangle mesh is challenging because scalar quantities like $\rho$ typically are discretized on vertices or faces while vectorial quantities like $J$ are better suited for edges.  Evaluating $\nicefrac{\|J\|^2}{2\rho}$ then requires averaging $J$ or $\rho$ so that the two end up on the same simplices.  If this problem is overcome, it still remains to prove a triangle inequality for discretizations of the Wasserstein distance resulting from such an approach.  Some recent papers with analogous constructions on graphs~\cite{maas2011gradient,solomon2016continuous,erbar2017computation} suggest that such an approach may be possible.

While the Benamou--Brenier dynamical formulation of transport is the best known, it is worth noting that the Beckmann problem~\eqref{eq:beckmann} for the 1-Wasserstein distance $\W_1$ more readily admits discretization using the finite element method (FEM) while preserving a triangle inequality.  Details of such a formulation as well as an efficient optimization algorithm are provided in~\cite{solomon2014earth}.  The reason~\eqref{eq:beckmann} is easier to discretize is that the time-varying aspect of transport is lost in this formulation:  All that is needed is a single vector $J(x)$ per point $x$.  What makes this problem easy to discretize and optimize is its downfall application-wise:  Interpolation with respect to $\W_1$ between two densities $\mu_0$ and $\mu_1$ is given by the uninteresting solution $\mu_t=(1-t)\mu_0+t\mu_1$, which does not displace mass but rather ``teleports'' it from the source to the target.

Another PDE-based approach to optimal transport is worth noting and has strong connections to the theory of transport without connecting to fluid flow.  Recall the Monge formulation of optimal transport on $\R^n$ in equation~\eqref{eq:monge}, which seeks a map $\phi(x)$ that pushes forward one distribution function $\rho_0(x)$ onto another $\rho_1(x)$.  A famous result by Brenier~\cite{brenier1991polar} shows that $\phi$ can be written as the gradient of a convex potential $\Psi(x)$:  $\phi(x)=\nabla\Psi(x)$.  Using $H$ to denote the Hessian operator, this potential satisfies the Monge-Amp\`ere PDE
\begin{equation}
\det(H\Psi(x))\rho_1(\nabla\Psi(x))=\rho_0(x),
\end{equation}
a second-order nonlinear elliptic equation that is extremely challenging to solve in practice.  A few algorithms, e.g.~\cite{oliker1989numerical,loeper2005numerical,benamou2010two,froese2011convergent,benamou2014numerical}, tackle this nonlinear system head-on, discretizing the variables involved and solving for $\Psi$.

\subsection{Semidiscrete Transport}\label{sec:semidiscrete_tr}

\begin{figure}
\centering
\includegraphics[width=.2\linewidth]{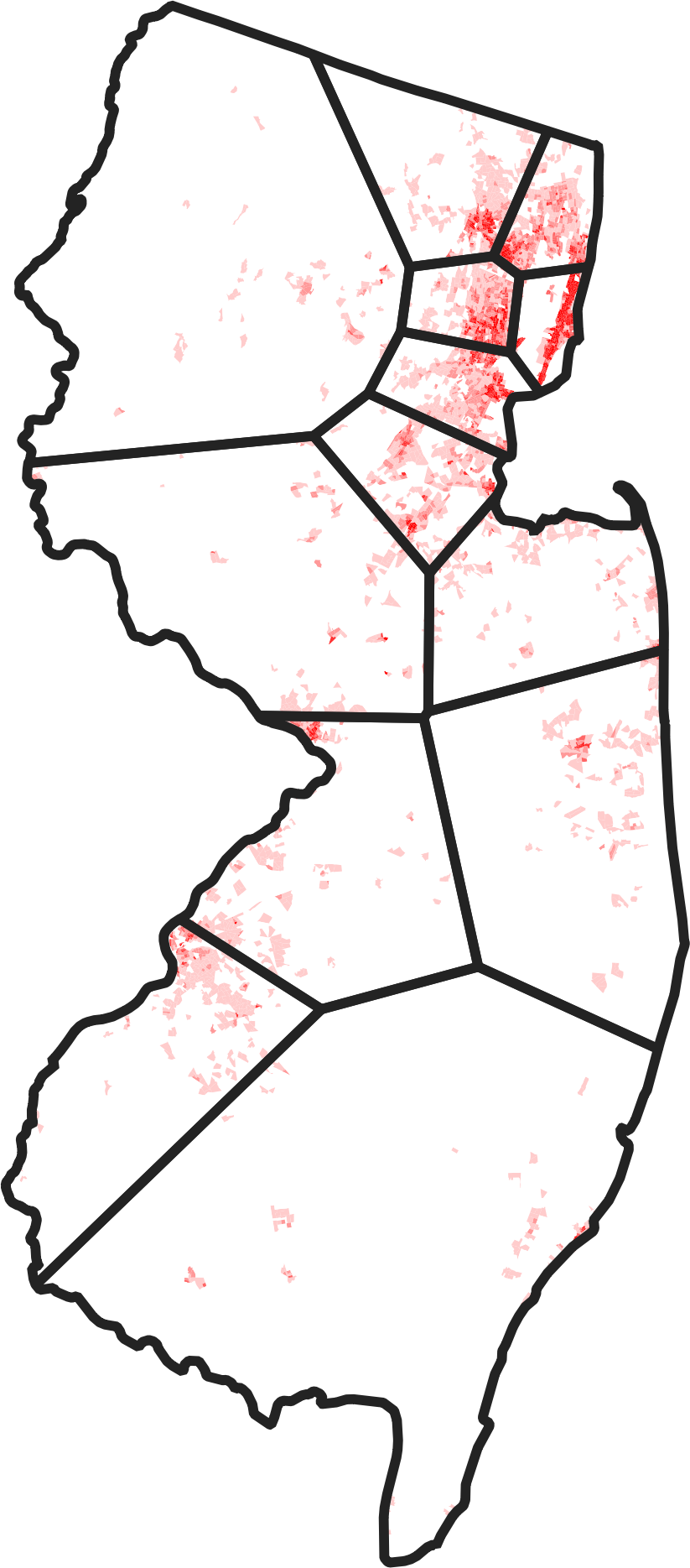}
\caption{Power diagram from a semidiscrete transport problem (image courtesy R.\ Barnes).  Here, semidiscrete transport is used to partition the state of New Jersey into cells with equal population; population density is shaded in red.}\label{fig:semidiscrete_example}
\end{figure}

Our final example from the computational transport world uses yet another formulation of the transport problem.  This time, our inspiration is the one-dimensional semidiscrete problem, whose solution is motivated from the formulation in equation~\eqref{eq:deltasum2}.  Our exposition of this problem closely follows the excellent tutorial~\cite{levy2017notions}.

In this setting, optimal transport is computed from a distribution whose mass is concentrated at a finite set of isolated points to a distribution with a known but potentially smooth density function.  Recall that in the one-dimensional case, we learned that each point of mass in the source is mapped to an \emph{interval} in the target.  That is, the domain of the target density  is partitioned into contiguous cells whose mass is assigned to a single source point.  We will find that the higher-dimensional analog is spiritually identical:  Each point of mass in the source density is assigned to a convex region of space  in the target.  This observation will suggest algorithms constructed from ideas in discrete geometry extending Voronoi diagrams and similar constructions.

As in~\eqref{eq:deltasum2}, suppose we are computing the 2-Wasserstein distance from a discrete measure $\mu:=\sum_{i=1}^ka_i\delta_{x_i}$, whose mass is concentrated at points $x_i\in\R^n$ with weights $a_i>0$, to an absolutely continuous measure $\nu$ with distribution function $\rho(x)$.  The dual formulation of transport~\eqref{eq:transportdual} in this case can be written
$$
\begin{array}{rl}
\sup_{\phi,\psi} & \sum_{i=1}^k a_i\phi(x_i) + \int_{\R^n} \psi(y)\rho(y)\,dA(y)\\
\st & \phi(x)+\psi(y)\leq c(x,y)\ \forall x,y\in\R^n.
\end{array}
$$
The objective in this case ``does not care'' about values of $\phi(x)$ for $x\not\in\{x_i\}_{i=1}^k.$  Define $\phi_i:=\phi(x_i).$  By this observation, we can write a problem with only one continuum variable:
$$
\begin{array}{rl}
\sup_{\phi,\psi} & \sum_i a_i\phi_i + \int_{\R^n} \psi(y)\rho(y)\,dA(y)\\
\st & \phi_i+\psi(y)\leq c(x_i,y)\ \forall y\in\R^n, i\in\{1,\ldots,k\}.
\end{array}
$$
In a slight abuse of notation, for the rest of this section we will think of $\phi$ as a vector $\phi\in\R^k$ rather than a function $\phi(x)$. Given the supremum, we might as well choose the largest $\psi$ possible that satisfies the constraints.  Hence,
$$
\psi(y)=\inf_{i\in\{1,\ldots,k\}} [c(x_i,y)-\phi_i].
$$
This leads to a final optimization problem in a \emph{finite} set of variables $\phi_1,\ldots,\phi_k$:
\begin{align}
\W_2^2(\mu,\nu)
&=
\sup_{\phi\in\R^k} \sum_i a_i\phi_i + \int_{\R^n} \rho(y)\left(\inf_{i\in\{1,\ldots,k\}} [c(x_i,y)-\phi_i]\right)\,dA(y)\nonumber\\
&=
\sup_{\phi\in\R^k} \sum_i \left[
a_i\phi_i + \int_{\Lag_\phi^c(x_i)} \rho(y)[c(x_i,y)-\phi_i] \,dA(y)
\right]\label{eq:semidiscrete_optim}
\end{align}
Here, $\Lag_\phi^c(x_i)$ indicates the \emph{Laguerre cell} corresponding to $x_i$:
\begin{equation}
\Lag_\phi^c(x_i):=\{y\in\R^n:
c(x_i,y)-\phi_i \leq c(x_j,y)-\phi_j\ \forall j\neq i
\}.
\end{equation}
The set of Laguerre cells yields the \emph{Laguerre diagram}, a partition of $\R^n$ determined by the cost function $c$ and the vector $\phi$; the $\phi_i$'s effectively control the sizes of the Laguerre cells in the diagram.  When $c(x,y)=\|x-y\|_2$ is a distance function and $\phi=0$, the Laguerre diagram equals the well-known Voronoi diagram of the $x_i$'s that partitions $\R^n$ into loci of points $S_i$ corresponding to those closer to $x_i$ than to the other $x_j$'s~\cite{aurenhammer1991voronoi}.  More importantly for the 2-Wasserstein distance, when $c(x,y)=\nicefrac{1}{2}\|x-y\|_2^2$, the Laguerre diagram is known as the \emph{power diagram}, an object studied since the early days of computational geometry~\cite{aurenhammer1987power}; an example is shown in Figure~\ref{fig:semidiscrete_example}.

Since~\eqref{eq:semidiscrete_optim} comes from a simplification of the dual of the transport problem, it is concave in $\phi$; a direct proof can be found in~\cite{aurenhammer1992minkowski}.  This implies that a simple gradient ascent procedure starting from any initial estimate of $\phi$ will reach a global optimum.  Define the objective function
$$F(\phi):=\sum_i\left[
a_i\phi_i + \int_{\Lag_\phi^c(x_i)} \rho(y)[c(x_i,y)-\phi_i] \,dA(y)
\right].$$
The gradient can be computed using the partial derivative expression
\begin{equation}\label{eq:semidiscrete_grad}
\frac{\partial F}{\partial \phi_i}=a_i-\int_{\Lag_\phi^c(x_i)}\rho(y)\,dA(y).
\end{equation}
This expression is predictable from the definition of $F(\phi)$; a similar formula exists for the second derivatives of $F$.  Setting the gradient~\eqref{eq:semidiscrete_grad} to zero formalizes 
 an intuition for the optimization problem~\eqref{eq:semidiscrete_optim}, that it resizes the Laguerre cells by modifying the $\phi_i$'s until the cell corresponding to each $x_i$ contains mass $a_i$:
$$
a_i=\int_{\Lag_\phi^c(x_i)}\rho(y)\,dA(y).
$$
%the intuition that $\phi_i$ is chosen so that the integral of $\phi(y)$ over the Laguerre cell corresponding to $x_i$ equals the mass $a_i$ assigned to $x_i$ in the discrete source distribution.

The main ingredient needed to compute the derivatives of $F$ is an algorithm for integrating $\rho$ over Laguerre cells.  Hence, gradient ascent and Newton's method applied to optimizing for $\phi$ cycle between updating the Laguerre diagram for the current $\phi$ estimate, recomputing the gradient and/or Hessian, assembling these into a search direction, and updating the current estimate of $\phi$.  For squared Euclidean costs, these algorithms are facilitated by fast algorithms for computing power diagrams, e.g.~\cite{bowyer1981computing,watson1981computing}.  While convergence of gradient descent with line search follows directly from concavity,~\cite{kitagawa2016newton} proves that under certain assumptions a damped version of Newton's algorithm---which employs the Hessian in addition to the gradient to accelerate convergence---exhibits global convergence.  

Example techniques following this template include~\cite{carlier2010knothe}, which proposed an early technique for 2D problems;~\cite{merigot2011multiscale}, for semi-discrete transport to piecewise-linear distribution functions in 2D supported on triangle meshes improved using a multiscale approximation; and~\cite{levy2015numerical}, which proposes semi-discrete transport to distributions in 3D that are piecewise-linear on tetrahedral meshes. \cite{de2012blue} provides an early example of a Newton solver for 2D semidiscrete transport using power diagrams and additionally uses derivatives of transport in the support points $x_i$ and weights $a_i$ for assorted approximation problems.  

Beyond providing fast algorithms for transport in the semidiscrete case, this formulation is also valuable for applications incorporating transport terms. 
\cite{de2011optimal} employs semidiscrete transport to a collection of distributions concentrated on line segments to reconstruct line drawings from point samples; \cite{digne2014feature} proposes a similar technique for reconstructing triangulated surfaces from point clouds in $\R^3$. 
\cite{goes2014weighted} defines a version of semi-discrete transport intrinsic to a triangulated surface, which can be used for tasks like parameterizing the set of per-vertex area weights in terms of the values $\phi_i$.

%\subsection{Randomized Transport}

%\subsection{Dualized Transport}

\section{Beyond Transport}

Beyond improving tools for solving the basic optimal transport problem, some of the most exciting recent work in computational transport involves using transport as a single term in a larger model.  In a recent tutorial for the machine learning community, we termed this new trend \emph{``Wassersteinization''}~\cite{cuturi2017primer}:  using Wasserstein distances to improve geometric properties of variational models in statistics, learning, applied geometry, and other disciplines.  Further extending the scope of applied transport, variations of the basic problem have been proposed to apply OT to objects other than probability distributions.

While a complete survey of these creative new applications and extensions is far beyond the scope of this tutorial, we highlight a few interesting pointers into the literature:
\begin{itemize}[leftmargin=*]
\item \textsc{Unbalanced transport:}  One limitation of the basic model for optimal transport is that it is a distance between histograms or probability distributions, rather than a distance between functions or vectors in $\R^n$---which may not integrate to $1$ or may contain negative values.  This leads to the problem of \emph{unbalanced transport}, in which mass conservation and/or positivity must be relaxed.  Models for this problem range from augmenting the transport problem with a ``trash can'' that can add or remove mass from distributions~\cite{pele2009fast} to extensions of dynamical transport to this case~\cite{chizat2016interpolating}.  Making transport work for functions rather than distributions while preserving the triangle inequality and other basic properties is challenging both theoretically and from a numerical perspective.
\item \textsc{Barycenters:}  The idea of displacement interpolation we motivated using~\eqref{eq:bb} suggests a generalization to more than two distributions, known as the \emph{Wasserstein barycenter} problem~\cite{agueh2011barycenters}.  Given $k$ distributions $\mu_1,\ldots,\mu_k$, the Wasserstein barycenter $\mu_{\mathrm{barycenter}}$ is defined as the minimizer of the following optimization problem 
\begin{equation}\label{eq:barycenter}
\mu_{\mathrm{barycenter}}:=\arg\min_{\mu}\sum_{i=1}^k \W_2^2(\mu,\mu_i).
\end{equation}
The Wasserstein barycenter gives some notion of \emph{averaging} a set of probability distributions, motivated by the observation that the average $\frac{1}{k}\sum_{i=1}^k x_i$ of a set of vectors $x_i\in\R^n$ is the minimizer $\arg\min_x\sum_i \|x-x_i\|_2^2$.  Barycenter algorithms range from extensions of the Sinkhorn algorithm~\cite{benamou2015iterative,solomon2015convolutional} to methods that perform gradient descent on $\mu$ by differentiating the distance $\W_2$ in its argument~\cite{cuturi2014fast} and stochastic techniques requiring only samples from the distributions $\mu_i$~\cite{staib2017parallel,claici2018stochastic}.  Other algorithms are inspired by a connection to multi-marginal transport~\cite{pass2015multi}, a generalization of optimal transport involving a distribution over the product of more than two measures.  The optimization problem~\eqref{eq:barycenter} is also one of the earliest examples of ``Wassersteinization,'' in the sense that it is an optimization problem for an unknown distribution $\mu$ including Wasserstein distance terms, contrasting somewhat from the optimization problems we considered in \S\ref{sec:discrete} in which the unknown is the transport distance itself. 

Further generalizing the barycenter problem leads to a notion of the Dirichlet energy of a map from points in one space to distributions over another~\cite{brenier2003extended,lavenant2017harmonic}, with applications in machine learning~\cite{solomon2014wasserstein} and shape matching~\cite{solomon2013dirichlet,mandad2017variance}. 
An intriguing recent paper also proposes an inverse problem for barycentric coordinates seeking weights for~\eqref{eq:barycenter} that ``explain'' an input distribution as a transport barycenter of others~\cite{bonneel2016wasserstein}.
\item \textsc{Quadratic assignment:}  The basic optimization problem for transport has an objective function that is \emph{linear} in the unknown transport matrix, expressing a preference for transport maps that do not move any single particle of probabilistic mass very far.  This model, however, does not necessarily extract \emph{smooth} maps, wherein distance traveled by any single particle is less important than making sure that nearby particles in the source are mapped to nearby locations in the target.  Such a smoothness term leads to a \emph{quadratic} term in the transport problem and allows it to be extended to a distance between metric-measure spaces known as the Gromov--Wasserstein distance~\cite{memoli2011gromov,memoli2014gromov}, inspired by the better-known but more rigid Gromov--Hausdorff distance.  From an optimization perspective, Gromov--Wasserstein computation leads to a ``quadratic assignment'' problem, known in the most general case to be NP-hard~\cite{sahni1976p}; practical instances of the problem in shape matching, however, can be tackled using spectral~\cite{memoli2009spectral} or entropy-based~\cite{solomon2016entropic} approximations and have shown promise for applications in shape matching.  \cite{peyre2016gromov} proposes a method for averaging metric spaces using a barycenter formulation similar to~\eqref{eq:barycenter}. 
\item \textsc{Capacity-constrained transport:}  Yet another extension of the transport problem comes from introducing \emph{capacity constraints} limiting the amount of mass that can travel between assorted pairs of source and target points; in the measure-theoretic formulation, this amounts to constraining transport plan to be dominated by another input plan~\cite{korman2015optimal}.  This constraint makes sense in many operations-type applications and has intriguing theoretical properties, but design of algorithms and discretizations for capacity-constrained transport remains largely open although~\cite{benamou2000computational} provides one approach again extending Sinkhorn's algorithm.
\item \textsc{Gradient flows and PDE:}  Given a function $f:M\rightarrow\R$ defined over a geometric space $M$ like a manifold, a \emph{gradient flow} of $f$ starting at some $x_0\in M$ attempts to minimize $f$ via ``gradient descent'' from $x(0):=x_0$ expressed as an ordinary differential equation (ODE) $x'(t)=-\nabla f(x(t)).$  Since OT puts a geometry on the space of distributions $\Prob(\R^n)$ over $\R^n$, we can define an analogous procedure that flows probability distributions to reduce certain functionals~\cite{jordan1998variational,santambrogio2017euclidean}.  For instance, gradient flow on the entropy functional~\eqref{eq:entropy} in the Wasserstein metric leads to the heat diffusion equation $\nicefrac{\partial \rho}{\partial t}=-\Delta \rho$, where $\Delta$ is the Laplacian operator; that is, performing gradient descent on entropy in the Wasserstein metric is exactly the same as diffusing the initial probability distribution like an unevenly-heated metal plate.  Beyond giving a variational motivation for certain PDE, this mathematical idea inspired numerical methods for solving PDE that can be written as gradient flows~\cite{peyre2015entropic,benamou2016augmented}.  Recent work has even incorporated transport into numerical methods for PDE that cannot easily be written as gradient flows in Wasserstein space, such as those governing incompressible fluid flow~\cite{li2010optimal,mirebeau2015numerical,de2015power,merigot2016minimal}.  Gradient flow properties can also be leveraged as structure to be preserved in discrete models of transport; for instance,~\cite{maas2011gradient} proposes a model for dynamical optimal transport on a graph and checks that the gradient flow of entropy---now an ODE rather than a PDE---agrees with a discrete heat equation.
\item \textsc{Matrix fields and vector measures:}  \emph{Vector measures} generalize probability measures by replacing scalar-valued probability values $\mu(S)\in[0,1]$ with values in other cones $\mathcal C$.  For instance, a tensor-valued measure $\mu$ assigns measurable sets $S$ to $d\times d$ postive semidefinite matrices $\mu(S)\in \mathcal S_+^d$ while satisfying analogous axioms to those laid out for probability measures in \S\ref{sec:transport_problem}.  These tensor fields find application in diffusion tensor imaging (DTI), which measures diffusivity of molecules like water in the interior of the human brain as a proxy for directionality of white matter fibers; OT extended to this setting can be used to align multiple such images.  A few recent models extend OT to this case and propose related numerical methods~\cite{ning2015matrix,chen2017matrix,peyre2017quantum}.
\end{itemize}

\section{Conclusion}

The techniques covered in this tutorial are just a few of many ways to approach discrete optimal transport.  New algorithms are proposed every month, and there is considerable room for mathematical, algorithmic, and application-oriented researchers to improve existing methods or make their own for different types of data or geometry.  Furthermore, mathematical properties such as convergence and approximation quality are still being established for new techniques.  Many questions also remain in linking to other branches of discrete differential geometry, e.g.\ at the most fundamental level defining a \emph{purely discrete} notion of optimal transport compatible with polyhedral meshes or simplicial complexes without requiring regularization and while preserving structure from the smooth case. 

These challenges aside, discrete optimal transport is demonstrating that OT holds interest far beyond mathematical analysis.  New discretizations and algorithms bring down OT's complexity to the point where it can be incorporated into practical engineering pipelines and into larger models without incurring a huge computational expense.  Further research into this new discipline holds unique potential to improve both theory and practice and eventually to bring insight into other branches of discrete and smooth geometry.

\subsection*{Acknowledgments.}  The author acknowledges the generous support of Army Research Office grant W911NF-12-R-0011 (``Smooth Modeling of Flows on Graphs''), from the MIT Research Support Committee (``Structured Optimization for Geometric Problems''), and from the MIT--IBM Watson AI Lab (``Large-Scale Optimal Transport for Machine Learning'').

Many thanks to MIT Geometric Data Processing Group members Mikhail Bessmeltsev, Edward Chien, Sebastian Claici, David Palmer, and Dima Smirnov for proofreading this document.

\bibliographystyle{amsplain}
\bibliography{optimal_transport}

\providecommand{\bysame}{\leavevmode\hbox to3em{\hrulefill}\thinspace}
\providecommand{\MR}{\relax\ifhmode\unskip\space\fi MR }
\providecommand{\MRhref}[2]{%
  \href{http://www.ams.org/mathscinet-getitem?mr=#1}{#2}
}
\providecommand{\href}[2]{#2}
\begin{thebibliography}{100}

\bibitem{Optimaldistricts}
\emph{{O}ptimal{D}istricts.org}, \url{http://www.optimaldistricts.org/}.

\bibitem{agueh2011barycenters}
Martial Agueh and Guillaume Carlier, \emph{Barycenters in the {W}asserstein
  space}, SIAM Journal on Mathematical Analysis \textbf{43} (2011), no.~2,
  904--924.

\bibitem{altschuler2017}
Jason Altschuler, Jonathan Weed, and Philippe Rigollet, \emph{Near-linear time
  approximation algorithms for optimal transport via {S}inkhorn iteration},
  Proc. NIPS, 2017, pp.~1961--1971.

\bibitem{arjovsky2017wasserstein}
Martin Arjovsky, Soumith Chintala, and L{\'e}on Bottou, \emph{Wasserstein
  generative adversarial networks}, International Conference on Machine
  Learning, 2017, pp.~214--223.

\bibitem{aurenhammer1987power}
Franz Aurenhammer, \emph{Power diagrams: properties, algorithms and
  applications}, SIAM Journal on Computing \textbf{16} (1987), no.~1, 78--96.

\bibitem{aurenhammer1991voronoi}
\bysame, \emph{Voronoi diagrams---a survey of a fundamental geometric data
  structure}, ACM Computing Surveys (CSUR) \textbf{23} (1991), no.~3, 345--405.

\bibitem{aurenhammer1992minkowski}
Franz Aurenhammer, Friedrich Hoffmann, and Boris Aronov, \emph{Minkowski-type
  theorems and least-squares partitioning}, Proceedings of the Eighth Annual
  Symposium on Computational Geometry, ACM, 1992, pp.~350--357.

\bibitem{bassetti2006minimum}
Federico Bassetti, Antonella Bodini, and Eugenio Regazzini, \emph{On minimum
  {K}antorovich distance estimators}, Statistics \& probability letters
  \textbf{76} (2006), no.~12, 1298--1302.

\bibitem{benamou2000computational}
Jean-David Benamou and Yann Brenier, \emph{A computational fluid mechanics
  solution to the {M}onge--{K}antorovich mass transfer problem}, Numerische
  Mathematik \textbf{84} (2000), no.~3, 375--393.

\bibitem{benamou2015iterative}
Jean-David Benamou, Guillaume Carlier, Marco Cuturi, Luca Nenna, and Gabriel
  Peyr{\'e}, \emph{Iterative {B}regman projections for regularized
  transportation problems}, SIAM Journal on Scientific Computing \textbf{37}
  (2015), no.~2, A1111--A1138.

\bibitem{benamou2016augmented}
Jean-David Benamou, Guillaume Carlier, and Maxime Laborde, \emph{An augmented
  {L}agrangian approach to {W}asserstein gradient flows and applications},
  ESAIM: Proceedings and Surveys \textbf{54} (2016), 1--17.

\bibitem{benamou2010two}
Jean-David Benamou, Brittany~D Froese, and Adam~M Oberman, \emph{Two numerical
  methods for the elliptic {M}onge-{A}mp\`ere equation}, ESAIM: Mathematical
  Modelling and Numerical Analysis \textbf{44} (2010), no.~4, 737--758.

\bibitem{benamou2014numerical}
\bysame, \emph{Numerical solution of the optimal transportation problem using
  the {M}onge--{A}mp\`ere equation}, Journal of Computational Physics
  \textbf{260} (2014), 107--126.

\bibitem{bernton2017inference}
Espen Bernton, Pierre~E Jacob, Mathieu Gerber, and Christian~P Robert,
  \emph{Inference in generative models using the {W}asserstein distance},
  arXiv:1701.05146 (2017).

\bibitem{bonneel2016wasserstein}
Nicolas Bonneel, Gabriel Peyr{\'e}, and Marco Cuturi, \emph{Wasserstein
  barycentric coordinates: histogram regression using optimal transport}, ACM
  Transactions on Graphics \textbf{35} (2016), no.~4, 71--1.

\bibitem{bowyer1981computing}
Adrian Bowyer, \emph{Computing {D}irichlet tessellations}, The Computer Journal
  \textbf{24} (1981), no.~2, 162--166.

\bibitem{boyd2011distributed}
Stephen Boyd, Neal Parikh, Eric Chu, Borja Peleato, and Jonathan Eckstein,
  \emph{Distributed optimization and statistical learning via the alternating
  direction method of multipliers}, Foundations and Trends in Machine Learning
  \textbf{3} (2011), no.~1, 1--122.

\bibitem{brenier1991polar}
Yann Brenier, \emph{Polar factorization and monotone rearrangement of
  vector-valued functions}, Communications on Pure and Applied Mathematics
  \textbf{44} (1991), no.~4, 375--417.

\bibitem{brenier2003extended}
\bysame, \emph{Extended {M}onge--{K}antorovich theory}, Lecture Notes in
  Mathematics (2003), 91--122.

\bibitem{carlier2010knothe}
Guillaume Carlier, Alfred Galichon, and Filippo Santambrogio, \emph{From
  {K}nothe's transport to {B}renier's map and a continuation method for optimal
  transport}, SIAM Journal on Mathematical Analysis \textbf{41} (2010), no.~6,
  2554--2576.

\bibitem{chen2017matrix}
Yongxin Chen, Tryphon~T Georgiou, and Allen Tannenbaum, \emph{Matrix optimal
  mass transport: a quantum mechanical approach}, IEEE Transactions on
  Automatic Control (2017).

\bibitem{chizat2016interpolating}
L{\'e}na{\"i}c Chizat, Gabriel Peyr{\'e}, Bernhard Schmitzer, and
  Fran{\c{c}}ois-Xavier Vialard, \emph{An interpolating distance between
  optimal transport and {F}isher--{R}ao metrics}, Foundations of Computational
  Mathematics (2016), 1--44.

\bibitem{claici2018stochastic}
Sebastian Claici, Edward Chien, and Justin Solomon, \emph{Stochastic
  {W}asserstein barycenters}, arXiv:1802.05757 (2018).

\bibitem{cohen1999earth}
Scott Cohen and Leonidas Guibas, \emph{The earth mover's distance under
  transformation sets}, Proc.\ ICCV, vol.~2, IEEE, 1999, pp.~1076--1083.

\bibitem{cohen2017balanced}
Vincent Cohen-Addad, Philip~N Klein, and Neal~E Young, \emph{Balanced power
  diagrams for redistricting}, arXiv:1710.03358 (2017).

\bibitem{courty2017optimal}
Nicolas Courty, R{\'e}mi Flamary, Devis Tuia, and Alain Rakotomamonjy,
  \emph{Optimal transport for domain adaptation}, PAMI \textbf{39} (2017),
  no.~9, 1853--1865.

\bibitem{crane2013geodesics}
Keenan Crane, Clarisse Weischedel, and Max Wardetzky, \emph{Geodesics in heat:
  A new approach to computing distance based on heat flow}, ACM Transactions on
  Graphics (TOG) \textbf{32} (2013), no.~5, 152.

\bibitem{cuturi2013sinkhorn}
Marco Cuturi, \emph{Sinkhorn distances: Lightspeed computation of optimal
  transport}, Advances in Neural Information Processing Systems, 2013,
  pp.~2292--2300.

\bibitem{cuturi2014ground}
Marco Cuturi and David Avis, \emph{Ground metric learning}, Journal of Machine
  Learning Research \textbf{15} (2014), no.~1, 533--564.

\bibitem{cuturi2014fast}
Marco Cuturi and Arnaud Doucet, \emph{Fast computation of {W}asserstein
  barycenters}, International Conference on Machine Learning, 2014,
  pp.~685--693.

\bibitem{cuturi2017primer}
Marco Cuturi and Justin Solomon, \emph{A primer on optimal transport}, NIPS
  Tutorial, 2017.

\bibitem{de2012blue}
Fernando De~Goes, Katherine Breeden, Victor Ostromoukhov, and Mathieu Desbrun,
  \emph{Blue noise through optimal transport}, ACM Transactions on Graphics
  (TOG) \textbf{31} (2012), no.~6, 171.

\bibitem{de2011optimal}
Fernando De~Goes, David Cohen-Steiner, Pierre Alliez, and Mathieu Desbrun,
  \emph{An optimal transport approach to robust reconstruction and
  simplification of 2d shapes}, Computer Graphics Forum, vol.~30, Wiley Online
  Library, 2011, pp.~1593--1602.

\bibitem{de2015power}
Fernando de~Goes, Corentin Wallez, Jin Huang, Dmitry Pavlov, and Mathieu
  Desbrun, \emph{Power particles: an incompressible fluid solver based on power
  diagrams}, ACM Transactions on Graphics \textbf{34} (2015), no.~4, 50--1.

\bibitem{digne2014feature}
Julie Digne, David Cohen-Steiner, Pierre Alliez, Fernando De~Goes, and Mathieu
  Desbrun, \emph{Feature-preserving surface reconstruction and simplification
  from defect-laden point sets}, Journal of Mathematical Imaging and Vision
  \textbf{48} (2014), no.~2, 369--382.

\bibitem{dobrushin1970definition}
Roland~L'vovich Dobrushin, \emph{Definition of random variables by conditional
  distributions}, Teoriya Veroyatnostei i ee Primeneniya \textbf{15} (1970),
  no.~3, 469--497.

\bibitem{douglas1956numerical}
Jim Douglas and Henry~H Rachford, \emph{On the numerical solution of heat
  conduction problems in two and three space variables}, Transactions of the
  American Mathematical Society \textbf{82} (1956), no.~2, 421--439.

\bibitem{erbar2017computation}
Matthias Erbar, Martin Rumpf, Bernhard Schmitzer, and Stefan Simon,
  \emph{Computation of optimal transport on discrete metric measure spaces},
  arXiv:1707.06859 (2017).

\bibitem{feydy2017optimal}
Jean Feydy, Benjamin Charlier, Fran{\c{c}}ois-Xavier Vialard, and Gabriel
  Peyr{\'e}, \emph{Optimal transport for diffeomorphic registration}, MICCAI
  2017, 2017.

\bibitem{ford1956solving}
Lester~Randolph Ford~Jr. and Delbert~Ray Fulkerson, \emph{Solving the
  transportation problem}, Management Science \textbf{3} (1956), no.~1, 24--32.

\bibitem{froese2011convergent}
Brittany~D Froese and Adam~M Oberman, \emph{Convergent finite difference
  solvers for viscosity solutions of the elliptic {M}onge--{A}mp{\`e}re
  equation in dimensions two and higher}, SIAM Journal on Numerical Analysis
  \textbf{49} (2011), no.~4, 1692--1714.

\bibitem{goes2014weighted}
Fernando~de Goes, Pooran Memari, Patrick Mullen, and Mathieu Desbrun,
  \emph{Weighted triangulations for geometry processing}, ACM Transactions on
  Graphics (TOG) \textbf{33} (2014), no.~3, 28.

\bibitem{haker2004optimal}
Steven Haker, Lei Zhu, Allen Tannenbaum, and Sigurd Angenent, \emph{Optimal
  mass transport for registration and warping}, International Journal of
  Computer Vision \textbf{60} (2004), no.~3, 225--240.

\bibitem{hirani2003discrete}
Anil~Nirmal Hirani, \emph{Discrete exterior calculus}, Ph.D. thesis, California
  Institute of Technology, 2003.

\bibitem{hitchcock1941distribution}
Frank~L Hitchcock, \emph{The distribution of a product from several sources to
  numerous localities}, Studies in Applied Mathematics \textbf{20} (1941),
  no.~1-4, 224--230.

\bibitem{jordan1998variational}
Richard Jordan, David Kinderlehrer, and Felix Otto, \emph{The variational
  formulation of the {F}okker--{P}lanck equation}, SIAM Journal on Mathematical
  Analysis \textbf{29} (1998), no.~1, 1--17.

\bibitem{kantorovich1942translocation}
Leonid~Vitalievich Kantorovich, \emph{On the translocation of masses}, Dokl.
  Akad. Nauk SSSR, vol.~37, 1942, pp.~199--201.

\bibitem{kitagawa2016newton}
Jun Kitagawa, Quentin M{\'e}rigot, and Boris Thibert, \emph{A {N}ewton
  algorithm for semi-discrete optimal transport}, arXiv:1603.05579 (2016).

\bibitem{klein1967primal}
Morton Klein, \emph{A primal method for minimal cost flows with applications to
  the assignment and transportation problems}, Management Science \textbf{14}
  (1967), no.~3, 205--220.

\bibitem{koopmans1941exchange}
Tjalling~C Koopmans, \emph{Exchange ratios between cargoes on various routes},
  (1941).

\bibitem{korman2015optimal}
Jonathan Korman and Robert McCann, \emph{Optimal transportation with capacity
  constraints}, Transactions of the American Mathematical Society \textbf{367}
  (2015), no.~3, 1501--1521.

\bibitem{kullback1951information}
Solomon Kullback and Richard~A Leibler, \emph{On information and sufficiency},
  The Annals of Mathematical Statistics \textbf{22} (1951), no.~1, 79--86.

\bibitem{kusner2015word}
Matt Kusner, Yu~Sun, Nicholas Kolkin, and Kilian Weinberger, \emph{From word
  embeddings to document distances}, International Conference on Machine
  Learning, 2015, pp.~957--966.

\bibitem{lavenant2017harmonic}
Hugo Lavenant, \emph{Harmonic mappings valued in the {W}asserstein space},
  arXiv:1712.07528 (2017).

\bibitem{levina2001earth}
Elizaveta Levina and Peter Bickel, \emph{The earth mover's distance is the
  {M}allows distance: Some insights from statistics}, Proc. ICCV, vol.~2, IEEE,
  2001, pp.~251--256.

\bibitem{levy2015numerical}
Bruno L{\'e}vy, \emph{A numerical algorithm for {$L_2$} semi-discrete optimal
  transport in 3{D}}, ESAIM: Mathematical Modelling and Numerical Analysis
  \textbf{49} (2015), no.~6, 1693--1715.

\bibitem{levy2017notions}
Bruno L{\'e}vy and Erica Schwindt, \emph{Notions of optimal transport theory
  and how to implement them on a computer}, Computers and Graphics \textbf{72}
  (2018), 135--148.

\bibitem{li2010optimal}
Bo~Li, Feras Habbal, and Michael Ortiz, \emph{Optimal transportation meshfree
  approximation schemes for fluid and plastic flows}, International Journal for
  Numerical Methods in Engineering \textbf{83} (2010), no.~12, 1541--1579.

\bibitem{lions1979splitting}
Pierre-Louis Lions and Bertrand Mercier, \emph{Splitting algorithms for the sum
  of two nonlinear operators}, SIAM Journal on Numerical Analysis \textbf{16}
  (1979), no.~6, 964--979.

\bibitem{loeper2005numerical}
Gr{\'e}goire Loeper and Francesca Rapetti, \emph{Numerical solution of the
  {M}onge--{A}mp{\`e}re equation by a {N}ewton's algorithm}, Comptes Rendus
  Mathematique \textbf{340} (2005), no.~4, 319--324.

\bibitem{lott2008some}
John Lott, \emph{Some geometric calculations on {W}asserstein space},
  Communications in Mathematical Physics \textbf{277} (2008), no.~2, 423--437.

\bibitem{maas2011gradient}
Jan Maas, \emph{Gradient flows of the entropy for finite {M}arkov chains},
  Journal of Functional Analysis \textbf{261} (2011), no.~8, 2250--2292.

\bibitem{mandad2017variance}
Manish Mandad, David Cohen-Steiner, Leif Kobbelt, Pierre Alliez, and Mathieu
  Desbrun, \emph{Variance-minimizing transport plans for inter-surface
  mapping}, ACM Transactions on Graphics \textbf{36} (2017), 14.

\bibitem{mccann1997convexity}
Robert~J McCann, \emph{A convexity principle for interacting gases}, Advances
  in Mathematics \textbf{128} (1997), no.~1, 153--179.

\bibitem{mccann2001polar}
\bysame, \emph{Polar factorization of maps on {R}iemannian manifolds},
  Geometric and Functional Analysis \textbf{11} (2001), no.~3, 589--608.

\bibitem{mccann1994convexity}
Robert~John McCann, \emph{A convexity theory for interacting gases and
  equilibrium crystals}, Ph.D. thesis, Princeton University, 1994.

\bibitem{memoli2009spectral}
Facundo M{\'e}moli, \emph{Spectral {G}romov--{W}asserstein distances for shape
  matching}, Proc. ICCV Workshops, IEEE, 2009, pp.~256--263.

\bibitem{memoli2011gromov}
\bysame, \emph{Gromov--{W}asserstein distances and the metric approach to
  object matching}, Foundations of Computational Mathematics \textbf{11}
  (2011), no.~4, 417--487.

\bibitem{memoli2014gromov}
\bysame, \emph{The {G}romov--{W}asserstein distance: A brief overview}, Axioms
  \textbf{3} (2014), no.~3, 335--341.

\bibitem{merigot2011multiscale}
Quentin M{\'e}rigot, \emph{A multiscale approach to optimal transport},
  Computer Graphics Forum, vol.~30, Wiley Online Library, 2011, pp.~1583--1592.

\bibitem{merigot2016minimal}
Quentin M{\'e}rigot and Jean-Marie Mirebeau, \emph{Minimal geodesics along
  volume-preserving maps, through semidiscrete optimal transport}, SIAM Journal
  on Numerical Analysis \textbf{54} (2016), no.~6, 3465--3492.

\bibitem{mikolov2013efficient}
Tomas Mikolov, Kai Chen, Greg Corrado, and Jeffrey Dean, \emph{Efficient
  estimation of word representations in vector space}, arXiv:1301.3781 (2013).

\bibitem{miller2007problem}
Stacy Miller, \emph{The problem of redistricting: the use of centroidal
  {V}oronoi diagrams to build unbiased congressional districts}, Senior
  project, Whitman College (2007).

\bibitem{mirebeau2015numerical}
Jean-Marie Mirebeau, \emph{Numerical resolution of {E}uler equations, through
  semi-discrete optimal transport}, Journ{\'e}es {\'E}quations aux
  D{\'e}riv{\'e}es Partielles (2015), 1--16.

\bibitem{mitchell1987discrete}
Joseph~SB Mitchell, David~M Mount, and Christos~H Papadimitriou, \emph{The
  discrete geodesic problem}, SIAM Journal on Computing \textbf{16} (1987),
  no.~4, 647--668.

\bibitem{monge1781memoire}
Gaspard Monge, \emph{M{\'e}moire sur la th{\'e}orie des d{\'e}blais et des
  remblais}, Histoire de l'Acad{\'e}mie Royale des Sciences de Paris (1781).

\bibitem{montavon2016wasserstein}
Gr{\'e}goire Montavon, Klaus-Robert M{\"u}ller, and Marco Cuturi,
  \emph{Wasserstein training of restricted {B}oltzmann machines}, Advances in
  Neural Information Processing Systems, 2016, pp.~3718--3726.

\bibitem{ning2015matrix}
Lipeng Ning, Tryphon~T Georgiou, and Allen Tannenbaum, \emph{On matrix-valued
  {M}onge--{K}antorovich optimal mass transport}, IEEE Transactions on
  Automatic Control \textbf{60} (2015), no.~2, 373--382.

\bibitem{oliker1987near}
Vladimir~I. Oliker, \emph{Near radially symmetric solutions of an inverse
  problem in geometric optics}, Inverse Problems \textbf{3} (1987), no.~4, 743.

\bibitem{oliker1989numerical}
Vladimir~I. Oliker and Laird~D. Prussner, \emph{On the numerical solution of
  the equation $\frac{{\partial ^2 z}}{{\partial x^2 }}\frac{{\partial ^2
  z}}{{\partial y^2 }} - \left( {\frac{{\partial ^2 z}}{{\partial x\partial
  y}}} \right)^2 = f$ and its discretizations, {I}}, Numerische Mathematik
  \textbf{54} (1989), no.~3, 271--293.

\bibitem{orlin1997polynomial}
James~B Orlin, \emph{A polynomial time primal network simplex algorithm for
  minimum cost flows}, Mathematical Programming \textbf{78} (1997), no.~2,
  109--129.

\bibitem{otto2001geometry}
Felix Otto, \emph{The geometry of dissipative evolution equations: the porous
  medium equation},  (2001).

\bibitem{papadakis2014optimal}
Nicolas Papadakis, Gabriel Peyr{\'e}, and Edouard Oudet, \emph{Optimal
  transport with proximal splitting}, SIAM Journal on Imaging Sciences
  \textbf{7} (2014), no.~1, 212--238.

\bibitem{pass2015multi}
Brendan Pass, \emph{Multi-marginal optimal transport: theory and applications},
  ESAIM: Mathematical Modelling and Numerical Analysis \textbf{49} (2015),
  no.~6, 1771--1790.

\bibitem{pele2009fast}
Ofir Pele and Michael Werman, \emph{Fast and robust earth mover's distances},
  Proc. ICCV, IEEE, 2009, pp.~460--467.

\bibitem{petersen2008matrix}
Kaare~Brandt Petersen and Michael~Syskind Pedersen, \emph{The matrix cookbook},
  Technical University of Denmark \textbf{7} (2008), 15.

\bibitem{peyre2010optimal}
Gabriel Peyr\'e, \emph{Optimal transport with {B}enamou--{B}renier algorithm},
  \url{http://www.numerical-tours.com/matlab/optimaltransp_2_benamou_brenier/},
  2010.

\bibitem{peyre2015entropic}
Gabriel Peyr{\'e}, \emph{Entropic approximation of {W}asserstein gradient
  flows}, SIAM Journal on Imaging Sciences \textbf{8} (2015), no.~4,
  2323--2351.

\bibitem{peyre2017quantum}
Gabriel Peyr\'e, L{\'e}na{\"i}c Chizat, Fran{\c{c}}ois-Xavier Vialard, and
  Justin Solomon, \emph{Quantum entropic regularization of matrix-valued
  optimal transport}, European Journal of Applied Mathematics (2017), 1--24.

\bibitem{peyre2017computational}
Gabriel Peyr{\'e} and Marco Cuturi, \emph{Computational optimal transport},
  Submitted, 2017.

\bibitem{peyre2016gromov}
Gabriel Peyr{\'e}, Marco Cuturi, and Justin Solomon,
  \emph{Gromov--{W}asserstein averaging of kernel and distance matrices},
  International Conference on Machine Learning, 2016, pp.~2664--2672.

\bibitem{plakhov2012billiards}
Alexander Plakhov, \emph{Billiards, optimal mass transport and problems of
  optimal aerodynamic resistance}, Preprint (2012).

\bibitem{ravindra1993network}
K.~Ahuja Ravindra, Thomas~L Magnanti, and James~B. Orlin, \emph{Network flows:
  theory, algorithms, and applications}, 1993.

\bibitem{rubner2000earth}
Yossi Rubner, Carlo Tomasi, and Leonidas~J Guibas, \emph{The earth mover's
  distance as a metric for image retrieval}, International journal of computer
  vision \textbf{40} (2000), no.~2, 99--121.

\bibitem{sahni1976p}
Sartaj Sahni and Teofilo Gonzalez, \emph{P-complete approximation problems},
  Journal of the ACM (JACM) \textbf{23} (1976), no.~3, 555--565.

\bibitem{santambrogio2015optimal}
Filippo Santambrogio, \emph{Optimal transport for applied mathematicians},
  Springer, 2015.

\bibitem{santambrogio2017euclidean}
\bysame, \emph{$\{${E}uclidean, metric, and {W}asserstein$\}$ gradient flows:
  an overview}, Bulletin of Mathematical Sciences \textbf{7} (2017), no.~1,
  87--154.

\bibitem{schwartzburg2014high}
Yuliy Schwartzburg, Romain Testuz, Andrea Tagliasacchi, and Mark Pauly,
  \emph{High-contrast computational caustic design}, ACM Transactions on
  Graphics (TOG) \textbf{33} (2014), no.~4, 74.

\bibitem{sethian1999fast}
James~A Sethian, \emph{Fast marching methods}, SIAM review \textbf{41} (1999),
  no.~2, 199--235.

\bibitem{sherman2017generalized}
Jonah Sherman, \emph{Generalized preconditioning and undirected minimum-cost
  flow}, Proc. SODA, SIAM, 2017, pp.~772--780.

\bibitem{sinkhorn1967concerning}
Richard Sinkhorn and Paul Knopp, \emph{Concerning nonnegative matrices and
  doubly stochastic matrices}, Pacific Journal of Mathematics \textbf{21}
  (1967), no.~2, 343--348.

\bibitem{slater1950lagrange}
Morton Slater, \emph{Lagrange multipliers revisited}, Cowles Commission
  Discussion Paper (1950), no.~403, 1--13.

\bibitem{solomon2018computational}
Justin Solomon, \emph{Computational optimal transport}, Snapshots of Modern
  Mathematics from Oberwolfach (2017), no.~8, 1--15.

\bibitem{solomon2015convolutional}
Justin Solomon, Fernando De~Goes, Gabriel Peyr{\'e}, Marco Cuturi, Adrian
  Butscher, Andy Nguyen, Tao Du, and Leonidas Guibas, \emph{Convolutional
  {W}asserstein distances: Efficient optimal transportation on geometric
  domains}, ACM Transactions on Graphics (TOG) \textbf{34} (2015), no.~4, 66.

\bibitem{solomon2013dirichlet}
Justin Solomon, Leonidas Guibas, and Adrian Butscher, \emph{Dirichlet energy
  for analysis and synthesis of soft maps}, Computer Graphics Forum, vol.~32,
  Wiley Online Library, 2013, pp.~197--206.

\bibitem{solomon2016entropic}
Justin Solomon, Gabriel Peyr{\'e}, Vladimir~G Kim, and Suvrit Sra,
  \emph{Entropic metric alignment for correspondence problems}, ACM
  Transactions on Graphics (TOG) \textbf{35} (2016), no.~4, 72.

\bibitem{solomon2014earth}
Justin Solomon, Raif Rustamov, Leonidas Guibas, and Adrian Butscher,
  \emph{Earth mover's distances on discrete surfaces}, ACM Transactions on
  Graphics (TOG) \textbf{33} (2014), no.~4, 67.

\bibitem{solomon2014wasserstein}
\bysame, \emph{Wasserstein propagation for semi-supervised learning},
  International Conference on Machine Learning, 2014, pp.~306--314.

\bibitem{solomon2016continuous}
\bysame, \emph{Continuous-flow graph transportation distances},
  arXiv:1603.06927 (2016).

\bibitem{staib2017parallel}
Matthew Staib, Sebastian Claici, Justin~M Solomon, and Stefanie Jegelka,
  \emph{Parallel streaming {W}asserstein barycenters}, Advances in Neural
  Information Processing Systems, 2017, pp.~2644--2655.

\bibitem{svec2007applying}
Lukas Svec, Sam Burden, and Aaron Dilley, \emph{Applying {V}oronoi diagrams to
  the redistricting problem}, The UMAP Journal \textbf{28} (2007), no.~3,
  313--329.

\bibitem{uzawa196810}
Hirofumi Uzawa, \emph{Iterative methods for concave programming}, Studies in
  Linear and Non-Linear Programming \textbf{2} (1968), 154.

\bibitem{varadhan1967behavior}
Sathamangalam R.\~Srinivasa Varadhan, \emph{On the behavior of the fundamental
  solution of the heat equation with variable coefficients}, Communications on
  Pure and Applied Mathematics \textbf{20} (1967), no.~2, 431--455.

\bibitem{vaserstein1969markov}
Leonid~Nisonovich Vaser\v{s}te\u{i}n, \emph{Markov processes over denumerable
  products of spaces, describing large systems of automata}, Problemy Peredachi
  Informatsii \textbf{5} (1969), no.~3, 64--72.

\bibitem{villani2003topics}
C{\'e}dric Villani, \emph{Topics in optimal transportation}, no.~58, American
  Mathematical Soc., 2003.

\bibitem{villani2008optimal}
\bysame, \emph{Optimal transport: old and new}, vol. 338, Springer Science \&
  Business Media, 2008.

\bibitem{wang1996design}
Xu-Jia Wang, \emph{On the design of a reflector antenna}, Inverse problems
  \textbf{12} (1996), no.~3, 351.

\bibitem{watson1981computing}
David~F Watson, \emph{Computing the $n$-dimensional {D}elaunay tessellation
  with application to {V}oronoi polytopes}, The Computer Journal \textbf{24}
  (1981), no.~2, 167--172.

\end{thebibliography}

\end{document}